\documentclass[dvips]{arxbj}
\usepackage{graphicx}
\usepackage{mathrsfs}

\aid{0}
\volume{16}
\issue{3}
\pubyear{2010}
\firstpage{679}
\lastpage{704}
\doi{10.3150/09-BEJ233}

\makeatletter

\newcommand{\iid}{\stackrel{\mathrm{iid}}{\sim}}
\newcommand{\mref}[1]{(\ref{#1})}
\newcommand{\eqref}[1]{(\ref{#1})}

\newcommand{\X}{\mathbf{X}}
\newcommand{\x}{\mathbf{x}}
\newcommand{\p}{\mathbf{p}}
\newcommand{\Z}{\mathbf{Z}}
\newcommand{\K}{\mathbf{K}}
\newcommand{\m}{\mathbf{m}}
\newcommand{\Y}{\mathbf{Y}}
\renewcommand{\S}{\mathbf{S}}
\newcommand{\q}{\mathbf{q}}
\newcommand{\bbeta}{\bolds{\beta}}
\newcommand{\bmu}{\bolds{\mu}}
\newcommand{\btau}{\bolds{\tau}}
\newcommand{\ind}{\stackrel{\mathrm{ind}}{\sim}}
\newcommand{\ddr}{\mathrm{d}}
\newcommand{\ssN}{N}
\newcommand{\rmP}{\mathrm{P}}

\newcommand{\rr}{\mathbb{R}}
\renewcommand{\ll}{L}
\newcommand{\pp}{\mathcal{P}}

\newproclaim{definition}{Definition}
\newtheorem{lemma}{Lemma}
\newremark{remark}{Remark}
\newremark{alg}{Algorithm}

\makeatother

\begin{document}
\begin{frontmatter}

\title{Bayesian nonparametric estimation and consistency of
mixed multinomial logit~choice models} \runtitle{Bayesian MMNL models}

\begin{aug}
\author[a]{\fnms{Pierpaolo} \snm{De Blasi}\corref{}\thanksref{a}\ead[label=e1]{pierpaolo.deblasi@unito.it}},
\author[b]{\fnms{Lancelot F.} \snm{James}\thanksref{b}\ead[label=e2]{lancelot@ust.hk}}
\and
\author[c]{\fnms{John W.} \snm{Lau}\thanksref{c}\ead[label=e3]{john@maths.uwa.edu.au}}
\runauthor{P. De Blasi, L.F. James and J.W. Lau}
\address[a]{Dipartimento di
Statistica e Matematica Applicata and Collegio Carlo Alberto, Universit\`a
degli Studi di Torino, corso Unione Sovietica 218/bis, 10134 Torino,
Italy.\\ \printead{e1}}
\address[b]{Department of Information Systems, Business Statistics and
Operations Management, The Hong Kong University of Science and
Technology,
Clear Water Bay, Kowloon, Hong Kong. \printead{e2}}
\address[c]{School of Mathematics and Statistics, University of
Western Australia, 35 Stirling
Highway, Crawley 6009, Western Australia, Australia. \printead{e3}}
\end{aug}

\received{\smonth{7} \syear{2008}} \revised{\smonth{4} \syear{2009}}

%
\begin{abstract}
This paper develops nonparametric estimation for discrete choice models
based on the mixed multinomial logit (MMNL) model. It has been shown
that MMNL models encompass all discrete choice models derived under the
assumption of random utility maximization, subject to the
identification of an unknown distribution $G$. Noting the mixture model
description of the MMNL, we employ a Bayesian nonparametric approach,
using nonparametric priors on the unknown mixing distribution $G$, to
estimate choice probabilities. We provide an important theoretical
support for the use of the proposed methodology by investigating
consistency of the posterior distribution for a general nonparametric
prior on the mixing distribution. Consistency is defined according to
an $L_1$-type distance on the space of choice probabilities and is
achieved by extending to a regression model framework a recent approach
to strong consistency based on the summability of square roots of prior
probabilities. Moving to estimation, slightly different techniques for
non-panel and panel data models are discussed. For practical
implementation, we describe efficient and relatively easy-to-use
blocked Gibbs sampling procedures. These procedures are based on
approximations of the random probability measure by classes of finite
stick-breaking processes. A simulation study is also performed to
investigate the performance of the proposed methods.
\end{abstract}

%

\end{frontmatter}

\section{Introduction}\label{section:1}
Discrete choice models arise naturally in many fields of application,
including marketing and transportation science. Such choice models are
based on the neoclassical economic theory of random utility
maximization (RUM). Given a finite set of choices $\mathbf{C}=\{
1,\ldots,J\} $, it is assumed that each individual has a utility function
\[
U_{j}={\mathbf{x}^{\prime} _{j}}{\bbeta}+\varepsilon_{j}\qquad
\mbox{for }j\in\mathbf{C}.
\]
The values $\mathbf{x}=( \mathbf{x}_{1},\ldots,\mathbf{x}_{J}) $ are observed covariates,
where $ \mathbf{x}_{j}\in\rr^{d}$ denote the covariates associated with
each choice $\{j\}\in\mathbf{C} $, the coefficient $ \bbeta$ is an
unknown (preference) vector in $\rr^d$ and $( \varepsilon
_{1},\ldots,\varepsilon_{J})$ are random terms. Suppose that all
$U_{j}$ are distinct and that the individual makes a choice $\{j\} $ if
and only if $U_{j}>U_{l}$ $\forall l\neq j$. The introduction of the
random error terms $\varepsilon_{j}$ represents a departure
from classical economic utility models. The random errors account
for the discrepancy between the actual utility, which is known by
the chooser, and that which is deduced by the experimenter who
observes $\mathbf{x}$ and the choice made by the individual. Hence, the
deterministic statement of choice $\{j\}$ is replaced by the
probability of choosing $\{j\}$, that is, $\rmP\{
U_{j}>U_{l}\ \forall l\neq j\}$.
The analysis of such a model depends on the specifications of the
errors. McFadden (\citeyear{mcfadden74}) shows that the specification of independent
Gumbel error terms leads to the tractable multinomial logit (MNL)
model. This representation is written as
\[
\rmP( \{j\}\vert \bbeta,\mathbf{x})=
\frac{\exp\{ \mathbf{x}_{j}^{\prime}\bbeta\} }{
\sum_{l\in\mathbf{C} }\exp\{ \mathbf{x}_{l}^{\prime}\bbeta
\} } \qquad\mbox{for }j\in\mathbf{C}.
\]
The MNL possesses the property of independence from irrelevant
alternatives (IIA), which makes it inappropriate in many situations.
The probit and the generalized extreme value models, which do not
exhibit the IIA
property and are models derived from dependent error structures, have been
proposed as alternatives to the MNL. A
drawback of the aforementioned procedures is that they are not
robust against model misspecification.

The mixed multinomial logit (MMNL) model, first introduced by
Cardell and Dunbar (\citeyear{cardell80}), emerges as potentially the most
attractive model. The book by Train (\citeyear{train03}) includes a detailed
discussion of this model. The general MMNL choice probabilities are
defined by mixing an MNL model over a mixing distribution $G$. For a
set of covariates $\mathbf{x}$, the MMNL model is written as
%
\begin{equation}\label{eq:MMNL}
\rmP( \{j\}\vert G,\mathbf{x})=
\int_{\rr^{d}} \frac{\exp\{ \mathbf{x}_{j}^{\prime}\bbeta\} }{
\sum_{l\in\mathbf{C} }\exp\{ \mathbf{x}_{l}^{\prime}\bbeta\}
} G(\ddr\bbeta)\qquad\mbox{for }j\in\mathbf{C}.
\end{equation}
McFadden and Train (\citeyear{McfTrain00}) establish the important result that, in
theory, all RUM models can be captured by correct specification of $G$.
Thus, a robust approach amounts to being able to employ statistical
estimation methods based on a nonparametric assumption on $G$. However,
statistical techniques have only been developed for the case where $G$
is given a parametric form. The most popular model is when $G$ is
specified to be multivariate normal with unknown mean $\bmu$ and
covariance matrix $\btau$:
%
\begin{equation}\label{eq:GML}
\rmP(\{j\}\vert\bmu,\btau,\mathbf{x})=
\int_{\rr^d} \frac{\exp\{ \mathbf{x}_{j}^{\prime}\bbeta\} }{
\sum_{l\in\mathbf{C} }\exp\{ \mathbf{x}_{l}^{\prime}\bbeta\} }
\phi(\bbeta|\bmu,\btau)\,\ddr\bbeta\qquad\mbox{for }j\in\mathbf{C},
\end{equation}
where $\phi( {\bbeta}|\bmu,\btau)$ represents a multivariate normal
density with parameters $\bmu$ and $\btau$. We shall refer to this as
a Gaussian mixed logit (GML) model. Here, based on a sample of size
$n$, one estimates the choice probabilities by estimating $\bmu$ and
$\btau$. Applications and discussions are found in, among others, Bhat
(\citeyear{Bhat98}), Brownstone and Train (\citeyear{brown99}), Erdem (\citeyear{erdem96}), Srinivasan and
Mahmassani (\citeyear{sri}) and Walker, Ben-Akiva and Bolduc (\citeyear{walkerJ07}).
Additionally, Dub\'{e} \textit{et al.} (\citeyear{dube02}) provide a discussion
focused on applications to marketing. The GML model is popular since
it is flexible and relatively easy to estimate via simulated maximum
likelihood techniques or via Bayesian MCMC procedures. Other
choices\vadjust{\goodbreak}
for $G$ include the lognormal and uniform distributions. Train
(\citeyear{train03}) discusses the merits and possible drawbacks of Bayesian MCMC
procedures versus simulated maximum likelihood procedures for
various choices of $G$.
However, despite the attractive features of the GML, it does not
encompass all RUM models, hence, it is not robust against
misspecification.

In this article, we develop a nonparametric Bayesian method for the
estimation of the choice probabilities and we prove consistency of
the posterior distribution. The idea is to model the mixing
distribution $G$ via a random probability measure in order to fully
exploit the flexibility of the MMNL model. Many nonparametric priors
are currently available for modeling $G$, such as stick-breaking priors,
normalized random measures with independent increments and Dirichlet
process mixtures. We establish consistency of the posterior
distribution of $G$ under neat sufficient conditions which are
readily verifiable for all of these nonparametric priors. Consistency
is defined according to an $\ll_1$-type distance on the space of
choice probabilities by exploiting the square root approach to
strong consistency of Walker (\citeyear{walk03a}, \citeyear{walk04}). We essentially show that
the Bayesian MMNL model is consistent if the prior on $G$ has the true
mixing distribution in its weak support and satisfies a mild condition on
the tails of the prior predictive distribution. We then move to
estimation and divide our discussion into methods for non-panel and
panel data. Specifically, for non-panel data models, we use, as a
prior for $G$, a mixture of Dirichlet processes. Methods for panel
data instead involve a Dirichlet mixture of normal densities. For
practical implementation, we describe efficient and relatively
easy-to-use blocked Gibbs sampling procedures, developed in Ishwaran and
Zarepour (\citeyear{IshZap00}) and Ishwaran and James
(\citeyear{IshJames01}).

The rest of the paper is organized as follows. In Section
\ref{section:2}, we describe the Bayesian nonparametric approach by
placing a nonparametric prior on the mixing distribution and present
the consistency result for the posterior distribution of $G$. In
Section \ref{section:3}, we show how to implement a blocked Gibbs
sampling for drawing inference when a discrete nonparametric prior
is used. Section \ref{section:4} deals with panel data with similar
Bayesian nonparametric methods, where we define a class of priors
for $G$ that preserves the distinct nature of individual preferences
and specialize the blocked Gibbs sampler to this setting. In Section
\ref{section:5}, we provide an illustrative simulation study which
shows the flexibility and good performance of our procedures.
Finally, in Section \ref{section:6}, we provide a detailed proof of
consistency.


%
%
\section{Bayesian MMNL models}\label{section:2}
A Bayesian nonparametric MMNL model is specified by placing
a nonparametric prior on the mixing distribution $G$ in
\eqref{eq:MMNL}:
%
\begin{equation}\label{eq:BayMMNL}
\rmP( \{j\}\vert \tilde{G},\mathbf{x})=
\int_{\rr^{d}} \frac{\exp\{ \mathbf{x}_{j}^{\prime}\bbeta\} }{
\sum_{l\in\mathbf{C} }\exp\{ \mathbf{x}_{l}^{\prime}\bbeta\}
} \tilde{G}(\ddr\bbeta)\qquad\mbox{for }j\in\mathbf{C}.
\end{equation}
Here, $\tilde{G}$ denotes a random probability measure which takes
values over the space $\mathbb{P}$ of probability measures on
$\rr^d$, the former endowed with the weak topology. The
nonparametric distribution of $\tilde{G}$ is denoted by $\pp$. Model
\eqref{eq:BayMMNL} can be equivalently expressed in hierarchical
form as
%
\begin{eqnarray}\label{eq:hierarchy}
Y_{i}\vert\bbeta_{i} & {\ind}&\frac{\exp\{
\mathbf{x}_{iY_{i}}^{\prime}\bbeta_{i}\} }{
\sum_{j\in\mathbf{C} }\exp\{
\mathbf{x}_{ij}^{\prime} \bbeta_{i}\} }\qquad\mbox{for }
i=1,\ldots,n\mbox{ and }Y_{i}\in\mathbf{C} ,
\nonumber\\
\bbeta_{i}\vert \tilde{G} & {\iid}& \tilde{G}
\qquad\mbox{for }i=1,\ldots,n, \\
\tilde{G} & {\sim}& \pp\nonumber
\end{eqnarray}
with $\mathbf{x}_i=( \mathbf{x}_{i1},\ldots,\mathbf{x}_{iJ}) $ the covariates
and $Y_i$ the
choice observed for individual $i$.

One can choose $\tilde{G}$ to be a Dirichlet process (Ferguson
(\citeyear{ferguson73})), although there currently exist other nonparametric priors that
can be used, like stick-breaking priors (Ishwaran and James (\citeyear{IshJames01}))
and normalized random measure with independent increments (NRMI)
(Regazzini, Lijoi and Pr\"{u}nster (\citeyear{RegLijPru03})).
All of these priors select discrete distributions almost surely (a.s.),
whereas random probability measures whose support contains
continuous distributions can be obtained by using a Dirichlet
process mixture of densities, in the spirit of Lo (\citeyear{lo84}). An
important role in the sequel will be played by the prior predictive
distribution of $\tilde{G}$, denoted by $H$, which is an element of
$\mathbb{P}$ and is defined by
%
\begin{equation}\label{eq:priorpred}
H(B):=\mathrm{E}[\tilde{G}(B)]
\end{equation}
for all Borel sets $B$ of $\rr^d$, where $\mathrm{E}(\cdot)$ denotes
expectation. In the next section, we show that an essential condition
for consistency of the posterior distribution
is expressed in terms of $H$. This yields an easy-to-use criterion
for the choice of the prior for $\tilde{G}$ as $H$ is readily
obtained for all of the nonparametric priors listed above. Furthermore,
one can embed a parametric model, such as the GML, within the
nonparametric framework via a suitable specification of the
distribution $H$.


\subsection{Posterior consistency}
Bayesian consistency deals with the asymptotic behavior of posterior
distributions with respect to repeated sampling. The problem can be
set in general terms as follows: suppose the existence of a ``true''
unknown distribution $P_0$ that generates the data, then check
whether the posterior accumulates in suitably-defined neighborhoods
of $P_0$. There exist two main approaches to the study of strong
consistency, that is, consistency when the neighborhood of $P_0$ is
defined according to the Hellinger metric on the space of density
functions. One is based on the metric entropy of the parameter space
and was set forth in Barron, Schervish and Wasserman (\citeyear{barron99}) and
Ghosal, Ghosh and Ramamoorthi (\citeyear{ghosal99}). The second approach was
introduced by Walker (\citeyear{walk03a}, \citeyear{walk04}) and has more of a Bayesian flavor, in
the sense that it relies on the summability of square roots of prior
probabilities. For discussion, the reader is referred to Wasserman
(\citeyear{wass98}), Walker, Lijoi and Pr\"{u}nster (\citeyear{WalkLijPr05}) and Choudhuri, Ghosal
and Roy (\citeyear{chouduri05}). Strong consistency in mixture models for density
estimation is addressed by Ghosal, Ghosh and Ramamoorthi (\citeyear{ghosal99}) and
Lijoi, Pr\"{u}nster and Walker (\citeyear{LijPruWal05}), by using the metric entropy
approach and the square root approach, respectively. As for the
non-identically distributed case, we mention Choi and Schervish
(\citeyear{choi07}) and Ghosal and Roy (\citeyear{GhosalRoy06}), both of which follow the metric
entropy approach. The square root approach is adopted by Walker (\citeyear{walk03b})
for nonparametric regression models and by Ghosal and Tang (\citeyear{GhosalTang06}) for
the estimation of transition densities in the context of Markov processes.

We face the issue of consistency for the MMNL model
\eqref{eq:BayMMNL} by exploiting the square root approach of Walker
and its variation proposed in Lijoi, Pr\"{u}nster and Walker (\citeyear{LijPruWal05})
which makes use of metric entropy in an instrumental way. We assume
the existence of a $G_0\in\mathbb{P}$ such that the true
distribution of $Y$ given $\X=\mathbf{x}$ is given by
\[
P_0(\{j\}\vert\mathbf{x})=\int_{\rr^d}\frac{\exp(\mathbf{x}_{j}^{\prime
}\bbeta)}{
\sum_{l\in\mathbf{C}}\exp(\mathbf{x}_{l}^{\prime}\bbeta)} G_0(\ddr
\bbeta).
\]
The variables $\X_i$ are taken as independent draws from a common
distribution $M(\ddr\mathbf{x})$ which is supported on $\mathcal{X}
\subset\rr^{Jd}$. The distribution of an infinite sequence
$(Y_i,\X_i)_{i\geq1}$ will be then denoted by
$\rmP^\infty_{(G_0,M)}$. Finally, let $\pp_n$ denote the posterior
distribution of $\tilde{G}$ given $(Y_1,\X_1),\ldots,(Y_n,\X_n)$;
see also equation \eqref{eq:Pi(A)} in Section \ref{section:6}. In the
sequel, we take the covariate distribution $M$ to be a fixed quantity
so that the posterior distribution does not depend on the specific form
of $M$. Note, however, that the posterior evaluation is also not
affected when $M$ is considered as a parameter with an independent
prior since it is reasonable to assume that the choice probabilities
are unrelated to $M$.

We give conditions on $G_0$ and the prior predictive distribution of
$\tilde{G}$ such that the posterior distribution $\pp_n$ concentrates
all probability mass in neighborhoods of $G_0$ defined according to
strong consistency of choice probabilities. To this end, we look at the
vector of choice probabilities as a vector-valued function $\q
\dvtx\mathcal{X}\to\Delta$, where $\Delta$ is the $J$-dimensional
probability simplex. We define
%
\begin{equation}\label{eq:choice}
\q(\mathbf{x}; G)= [\rmP(\{1\}\vert G,\mathbf{x}),\ldots,\rmP(\{J\}
\vert G,\mathbf{x}) ]
\end{equation}
for any $G\in\mathbb{P}$. On the space $\mathscr{Q}=\{\q(\cdot; G)\dvt
G\in\mathbb{P}\}$, we define the $\ll_1$-type distance
%
\begin{equation}\label{eq:hellinger}
d(\q_1,\q_2)=
\int_{\mathcal{X}}|\q_1(\mathbf{x})-\q_2(\mathbf{x})|M(\ddr\mathbf{x}),
\end{equation}
where $|\cdot|$ denotes the Euclidean norm in $\Delta$.
%
%
\begin{definition}
$\pp$ is consistent at $G_0$ if, for any $\epsilon>0$,
\[
\pp_n\{G\dvt d (\q(\cdot; G),\q(\cdot; G_0) )>\epsilon\}\to0,\qquad
\rmP_{(G_0,M)}^{\infty}\mbox{-a.s.}
\]
\end{definition}
The main result is stated in the following theorem.
%
%
\begin{thm}\label{theorem:consistency}
Let $\pp$ be a prior on $\mathbb{P}$ with predictive distribution
$H$ and $G_0$ be in the weak support of $\pp$. Suppose that $\mathcal
{X}$ is a compact subset of $\rr^{Jd}$. If
\begin{longlist}[(ii)]
\item[(i)] $P_0(\{j\}\vert\mathbf{x})>0$ for any $j\in\mathbf{C}$
and $\mathbf{x}\in\mathcal{X}$;
\item[(ii)] $\int_{\rr^d}|\bbeta|H(\ddr\bbeta)<+\infty$,
\end{longlist}
then $\pp$ is consistent at $G_0$.
\end{thm}

The compactness of the covariate space is a standard assumption in
nonparametric regression problems. Condition (i) is fairly
reasonable since it is guaranteed by a correct specification of the
RUM model: one can always redefine the set of choices or the covariate
space to fulfill this requirement. Moreover, because of the compactness
of $\mathcal{X}$, condition (i) implies that $G_0$ is a proper
distribution on $\rr^d$, that is, with no masses escaping at infinity.
The verification that $G_0$ belongs to the weak support of $\pp$ is
then an easy task: in general, it is sufficient that the prior
predictive distribution $H$ has full support on $\rr^d$.
Condition (ii) is a mild condition on the tails of $H$: it is
satisfied by any distribution with tails lighter than the Cauchy distribution.


\subsection{Illustration}
It is worth considering condition (ii) in more detail for a
variety of Bayesian MMNL models, obtained from different
specifications of $\pp$.
%
%
If $\tilde{G}$ is taken to be a Dirichlet process with base
measure $\alpha=aF$, where $a>0$ is a constant and $F\in\mathbb{P}$,
then $F$ coincides with $H$ in \eqref{eq:priorpred}.
%
%
A larger class of Bayesian MMNL models arise when $\tilde{G}$ is
chosen to be a stick-breaking prior:
%
\begin{equation}\label{eq:stick}
\tilde{G}(\cdot)=\sum_{k\geq1}p_{k}\delta_{Z_{k}}(\cdot),
\end{equation}
where the $p_k$ are positive random probabilities chosen to be
independent of $Z_{k}$ and such that $\sum_{k\geq1} p_k=1$ a.s.
The $Z_k$ are random locations taken as independent draws from
some non-atomic distribution $F$ in $\mathbb{P}$. What characterizes
a stick-breaking prior is that the random weights are expressible as
$p_{k}=V_{k}\prod_{i=1}^{k-1}(1-V_{i})$, where the $V_{k}$ are
independent beta-distributed random variables of parameters
$a_k,b_k>0$; we write $V_k\sim\operatorname{beta}(a_k,b_k)$.
Examples of random probability measures in this class are given in
Ishwaran and James (\citeyear{IshJames01}); see also Pitman and Yor (\citeyear{PitYor97}) and
Ishwaran and Zarepour (\citeyear{IshZap00}). They represent extensions of the
Dirichlet process, which has $a_k=1$ and $b_k=a$ $\forall k$, and
they all have in common that the prior predictive distribution $H$
coincides with $F$.

%
%
The class of NRMI is another valid choice for $\pp$. Specifically,
one can take $\tilde{G}(\cdot)=\tilde\mu(\cdot)/\break\tilde\mu(\rr^d)$,
where $\tilde\mu$ is a completely random measure with Poisson
intensity measure $\nu(\ddr v,\ddr z)=\rho(\ddr v| z)\alpha(\ddr z)$
on $(0,+\infty)\times\rr^d$. Here, $\rho(\cdot| z)$ is a L\'{e}vy
density on $(0,+\infty)$ for any $z$ and $\alpha$ is a finite
measure on $\rr^d$ such that
$\psi(u):=\int_{\rr^d\times\rr^+}(1-\mathrm{e}^{-uv})\rho(\ddr v|
z)\alpha(\ddr z)<\infty$, which is needed to guarantee that
$\tilde\mu(\rr^d)<\infty$ a.s. It can be shown that
$H(B)=\int_B\int_0^{+\infty}\mathrm{e}^{-\psi(u)} \{
\int_0^{+\infty}\mathrm{e}^{-uv}\times v\rho(\ddr v| z) \}
\,\ddr u \,\alpha(\ddr z)$
for any Borel set $B$ of $\rr^d$; see also James, Lijoi and
Pr\"{u}nster (\citeyear{JamLijPru05}).
When $\rho(\ddr v| z)=\rho(\ddr v)$ for each $z$ (homogeneous case),
the prior predictive distribution reduces to
%
\begin{equation}\label{eq:species}
H(B)=\frac{\alpha(B)}{\alpha(\rr^d)}\qquad\mbox{for any Borel
}B\subset\rr^d.
\end{equation}
The homogeneous NRMI includes, as a special case, the Dirichlet
process and belongs, together with the stick-breaking priors, to the
class of \textit{species sampling models}, for which \eqref{eq:species}
holds for some finite measure $\alpha$.
Note that all of the nonparametric priors belonging to this class allow
an easy verification of condition (ii).

%
%
The specification of the nonparametric prior in terms of a base
measure $\alpha$, as in \eqref{eq:species}, allows more
flexibility to be introduced via an additional level in the hierarchal structure
\eqref{eq:hierarchy}. If we let the base measure be indexed by a
parameter $\theta$, say $\alpha_\theta$, and $\theta$ be random with
probability density $\pi(\theta)$ on some Euclidean space $\Theta$,
then we obtain a mixture of Dirichlet process in the spirit of
Antoniak (\citeyear{anto74}).
Condition (ii) must then be verified for the convolution
%
\begin{equation}\label{eq:Antoniak}
H(B)=\int_\Theta\int_B H_\theta(\ddr z)
\pi(\theta)\,\ddr\theta,\qquad\mbox{where }
H_\theta(\ddr z)=\frac{\alpha_\theta(\ddr z)}{\alpha_\theta(\rr^d)}.
\end{equation}
It is quite straightforward to check that condition (ii) holds for
the mixture of Dirichlet processes implemented in the analysis of
non-panel data in Section \ref{section:3}.

%
%
Finally, consider the case of Dirichlet process mixture models of Lo
(\citeyear{lo84}), where $\tilde{G}$ is absolutely continuous with respect to
the Lebesgue measure on $\rr^{d}$
with random density function specified as
$\int_\Theta K(\bbeta,\theta)\tilde{\Pi}(\ddr\theta)$. Here,
$K(\bbeta,\theta)$ is a non-negative kernel defined on
$\rr^d\times\Theta$ such that, for each $\theta\in\Theta$,
$\int_{\rr^d}K(z,\theta)\,\ddr z=1$, while $\tilde{\Pi}$ is a
Dirichlet process prior with base measure $aF$ and $F$ a probability
measure on $\Theta$. The distribution $H$ is then absolutely
continuous and is given by
\[
H(B)=\int_B\int_\Theta K(z,\theta)F(\ddr\theta)\,\ddr z.
\]
As in \eqref{eq:Antoniak}, verifying condition (ii) requires a study
of the tail properties of a convolution, this time of $K(z,
\theta)$ with respect to $F(\ddr\theta)$. In the analysis of panel
data (see Section \ref{section:4}), we adopt a Dirichlet mixture
model as continuous nonparametric prior for $\tilde{G}$ where the
verification of condition (ii) can be readily established.


%
%
\section{Implementation for non-panel data}\label{section:3}
Assume that we have a single observation for each individual and that
we want
to account for the possibility of ties among different individuals'
preferences. Therefore, we use a discrete nonparametric prior for
the mixing distribution.
Take $\tilde{G}$ to be a Dirichlet process with base measure $a F$
and denote its law by $\mathcal{P}(\ddr G| a F)$ (although the
treatment can be easily extended to any other stick-breaking prior).
Representation \eqref{eq:stick} then holds with random
probabilities $p_1,p_2,\ldots$ at locations $Z_1,Z_2,\ldots,$ which
are i.i.d.~draws from $F$. This translates into a Bayesian model for
the MMNL as
%
\begin{equation}\label{eq:BMMNL}
\rmP(\{j\}\vert\tilde{G},\mathbf{x}) =
\sum_{k\geq1}p_k\frac{ \exp\{\mathbf{x}_j^{\prime}Z_k\} }{
\sum_{l\in\mathbf{C}}\exp\{\mathbf{x}_l^{\prime}Z_k\} }\qquad\mbox{for
}j\in\mathbf{C}.
\end{equation}
One can then center $\tilde{G}$ on a parametric model like the GML
in \eqref{eq:GML} by taking $F$ to have normal density
$\phi(\bbeta|\bmu,\btau)$.
In a parametric Bayesian framework, by placing priors on
$\bmu,\btau$, one is able to get posterior estimates of
$\bmu,\btau$, but inference is restricted to the assumption of the
GML model. The flexibility of the Bayesian nonparametric approach
allows one to choose $F$ based on convenience and ease of use and to
utilize, for instance, the attractive features of GML models while
still maintaining the robustness of a nonparametric approach.

In the case of the Dirichlet process, the parameters associated with
$F$, for instance, $\bmu$ and $\btau$, are considered fixed. As
observed in Section \ref{section:2}, one can introduce more
flexibility in the model by treating such parameters as random.
Specifying $\theta=(\bmu,\btau)$, $F_\theta(\ddr\bbeta)$ to have
density $\phi(\bbeta|\theta)\,\ddr\bbeta$ and $\pi(\theta)$ to be the
density function for $\theta$, the law of $\tilde{G}$ is given by
the mixture $\int_{\Theta} {\mathcal P}(\ddr G|a
F_{\theta})\pi(\ddr\theta)$. Equivalently, using \mref{eq:stick}, a
mixture of Dirichlet processes is defined by specifying each
$Z_{k}\vert\theta$ to be i.i.d.~$F_{\theta}$. Note that, conditional
on $\theta$, a prior guess for the choice probabilities is
%
\begin{equation}\label{eq:MMNLtheta}
\mathrm{E} [\rmP(\{j\}\vert \tilde{G},\mathbf{x})\vert\theta]=
\int_{\rr^{d}}\frac{\exp\{\mathbf{x}_j^{\prime}\bbeta\}}{\sum_{l\in
\mathbf{C}}
\exp\{\mathbf{x}_l^{\prime}\bbeta\} }F_\theta(\ddr\bbeta)\qquad\mbox{for
}j\in
\mathbf{C}.
\end{equation}
By the properties of the Dirichlet process, the prediction rule for
the choice probabilities given $\bbeta_{1},\ldots, \bbeta_{n}$ is
given by
\begin{eqnarray}\label{eq:prediction}
&&\mathrm{E} [\rmP(\{j\}\vert \tilde{G},\mathbf{x})\vert\theta,\bbeta
_{1},\ldots, \bbeta_{n} ]\nonumber\\[-8pt]\\[-8pt]
&&\quad=\frac{a}{a+n}\rmP(\{j\}\vert F_{\theta},\mathbf{x})+\sum
_{i=1}^{n}\frac{1}{a
+n} \frac{\exp\{ \mathbf{x}_{j}^{\prime}\bbeta_{i}\} }{
\sum_{l\in\mathbf{C} }\exp\{ \mathbf{x}_{l}^{\prime}\bbeta_{i} \}},\nonumber
\end{eqnarray}
where $\rmP(\{j\}\vert F_{\theta},\mathbf{x}):=\mathrm{E} [\rmP(\{j\}
\vert
\tilde{G},\mathbf{x})\vert\theta]$ is given in \eqref{eq:MMNLtheta} with
a notation consistent with \eqref{eq:MMNL}. However, the variables
$\bbeta_{i}$ are not observable and hence one needs to implement
computational procedures to draw from their posterior distribution.

In this framework, a reasonable algorithm to use is the blocked
Gibbs sampler developed in Ishwaran and Zarepour (\citeyear{IshZap00}) and Ishwaran
and James (\citeyear{IshJames01}). Indeed, since the multinomial logistic kernel does
not form a conjugate pair for $\bbeta$, marginal algorithms suffer
from slow convergence, although strategies for overcoming this
problem can be found in MacEachern and Muller (\citeyear{MachMull98}).


\subsection{Blocked Gibbs algorithm}
In this section, we discuss how to implement a blocked Gibbs sampling
algorithm for drawing inference on a nonparametric hierarchical
model with the structure
\begin{eqnarray}\label{eq:model3}
Y_i\vert\bbeta_i & \ind&L(Y_i,\bbeta_i)\qquad
\mbox{for }i=1,\ldots,n\mbox{ and }Y_i\in\mathbf{C},\nonumber\\
\bbeta_{i}\vert\tilde{G} & \iid& \tilde{G}\qquad\mbox{for }i=1,\ldots
,n,\nonumber\\[-8pt]\\[-8pt]
\tilde{G}\vert\theta& \sim& \pp(\ddr G|a F_\theta),\nonumber\\
\theta&\sim&\pi(\ddr\theta),\nonumber
\end{eqnarray}
where $L(Y_i,\bbeta)=\exp\{\mathbf{x}_{iY_i}^\prime\bbeta\}/\sum
_{j\in
\mathbf{C}}\exp\{\mathbf{x}_{ij}^\prime\bbeta\}$ is the probability
for $Y_i$\vspace*{1pt}
conditional on $\bbeta_i$. The blocked Gibbs sampler utilizes the
fact that a truncated Dirichlet process, discussed in Ishwaran and
Zarepour (\citeyear{IshZap00}) and Ishwaran and James (\citeyear{IshJames01}), serves as a good
approximation to the random probability measure $\tilde G\vert\theta$
in \eqref{eq:model3}.
We replace the conditional law $\pp(\ddr G|a F_\theta)$ with the law
of the random probability measure
%
\begin{equation}\label{eq:dpn}
\tilde G(\cdot)=\sum_{k=1}^{N}p_{k}\delta_{Z_{k}}(\cdot),\qquad 1\leq N<\infty,
\end{equation}
where $Z_{k}\vert\theta$ are i.i.d.~$F_{\theta}$ and the random
probabilities $p_{1},\ldots,p_{N}$ are defined by the stick-breaking
construction
%
\begin{equation}\label{eq:finitestick}
p_1=V_1 \quad\mbox{and}\quad
p_k=(1-V_1)\cdots(1-V_{k-1})V_k,\qquad k=2,\ldots,N,
\end{equation}
with $V_1,V_2,\ldots,V_{N-1}$ i.i.d.~$\operatorname{beta}(1,a)$ and $V_{N}=1$,
which ensures that $\sum_{k=1}^{N}p_{k}=1$.
The law of $\tilde G\vert \theta$ in \eqref{eq:dpn} is referred to as
a \textit{truncated Dirichlet process} and will be denoted
$\pp^{\ssN}(\mathrm{d}G|\alpha F_{\theta})$. Moreover, the limit as
$N\to\infty$ will converge to a random probability measure with law
$\pp(\ddr G|a F_\theta)$. Indeed, the method yields an accurate
approximation of the Dirichlet process for $N$ moderately large
since the truncation is exponentially accurate.
Theorem 2 in Ishwaran and James (\citeyear{IshJames01}) provides an $\ll_{1}$-error
bound for the approximation of conditional density of
$\Y=(Y_1,\ldots,Y_n)$ given $\theta$. Let
\[
\mu^{\ssN}(\Y|\theta)=\int \Biggl[\prod_{i=1}^n\int_{\rr^d}
L(Y_i,\bbeta_i)
G(\ddr\bbeta_i) \Biggr] \pp^{\ssN}(\ddr G| a F_\theta)
\]
and $\mu(\Y|\theta)$ be its limit under the prior $\pp(\ddr G|a
F_\theta)$. One then has
\[
\| \mu^{\ssN}-\mu \| _{1}:=
\int\bigl|\mu^{\ssN}(\Y|\theta)-\mu(\Y|\theta) \bigr|
\,\ddr\Y\sim4n\mathrm{e}^{- ( N-1 ) / a },
\]
where the integral above is considered over the counting measure on
the $n$-fold product space $\mathbf{C}^n$.
Moreover, Corollary 1 in Ishwaran and James (\citeyear{IshJames02}) can be used to show
that the truncated Dirichlet process also leads to asymptotic
approximations to the posterior that are exponentially accurate.

The key to working with random probability measures like
\eqref{eq:dpn} is that it allows blocked updates to be performed for
$\p=(p_1,\ldots,p_n)$ and $\Z=(Z_1,\ldots,Z_n)$ by recasting the
hierarchical model \eqref{eq:model3} completely in terms of random
variables. To this aim, define the classification variables
$\K=\{K_{1},\ldots,K_{n}\}$ such that, conditional on $\p$, each
$K_{i}$ is independent with distribution
\\[-15pt]
\[
\rmP\{K_i\in\cdot\vert\p\}=\sum_{k=1}^{N}p_{k}\delta_{k}(\cdot).
\]
That is, $\rmP\{K_i=k\vert\p\}=p_{k}$ for $k=1,\ldots, N$ so that
$K_i$ identifies the $Z_k$ associated with each $\bbeta_i$:
$\bbeta_{i}=Z_{K_{i}}$. In this setting, a sample
$\bbeta_{1},\ldots,\bbeta_{n}$ from \eqref{eq:dpn} produces
$n_{0}\leq\min(n,N)$ distinct values. The blocked Gibbs
algorithm is based on sampling $\K,\p,\Z,\theta$ from the
distribution proportional to
\[
\Biggl[\prod_{i=1}^{n}L(Y_i,\bbeta_i) \Biggr]
\Biggl[\prod_{i=1}^{n}\sum_{k=1}^Np_{k}\delta_{Z_k}(\ddr\bbeta_i) \Biggr]
\pi(\p) \Biggl[\prod_{k=1}^{N}F_\theta(\ddr Z_{k}) \Biggr]\pi(\ddr\theta),
\]
where $\pi(\p)$ denotes the distribution of $\p$ defined in
\eqref{eq:finitestick}. This augmented likelihood is an expression
of the augmented density when $\pp(\ddr G|a F_{\theta})$ is
replaced by~$\pp^{\ssN}(\ddr G|a F_{\theta})$.

Before describing the algorithm, we specify choices for $F_{\theta}$
and $\theta$ which agree with the GML model. Set
$\theta=(\bmu,\btau)$ and specify the density of $F_{\theta}$ to be
$\phi(\bbeta|\bmu,\btau)$. Let $\lambda$ denote a positive scalar.
We choose a multivariate normal inverse Wishart distribution for
$\bmu,\btau$, where, specifically, $\bmu\vert\btau$ is a multivariate
normal vector with mean parameter $\m$ and scaled covariance matrix
$\lambda^{-1}\btau$ and $\btau$ is drawn from an inverse Wishart
distribution with degrees of freedom $\nu_{0}$ and scale matrix
$\S_{0}$. We denote this distribution for $\bmu,\btau$ as
$\operatorname{N\mbox{-}IW}(\m,\lambda^{-1}\btau,\nu_0,\S_0)$. Our specification is
similar to that used in Train (\citeyear{train03}), Chapter 12, for a parametric
GML model for panel data.

\begin{alg}\label{alg1}
\begin{enumerate}
\item\textit{Conditional draw for $\K$.} Independently sample
$K_i$ according to $\rmP\{K_i\in\cdot|\p,\Z,\Y\}=\sum_{k=1}^N
p_{k,i}\delta_{k}(\cdot)$ for $i=1,\ldots,n$, where
\[
(p_{1,i},\ldots,p_{N,i})\propto (p_1 L(Y_i,Z_1),\ldots,p_N
L(Y_i,Z_N) ).
\]
\item\textit{Conditional draw for $\p$.} $p_1=V_1^*$,
$p_k=(1-V_1^*)\cdots(1-V_{k-1}^*)V_k^*$, $k=2,\ldots,N-1$ and
$V_N^*=1$, where, if $e_k$ records the number of $K_i$ values which
equal $k$,
\[
V_k^*\ind\operatorname{beta} \Biggl(1+e_k,a+\sum_{l=k+1}^{N}e_l
\Biggr),\qquad k=1,\ldots,N-1.
\]
\item\textit{Conditional draw for $\Z$.} Let
$\{K_1^*,\ldots,K_{n_0}^*\}$ denote the unique set of $K_i$
values.\\
For each $k\notin\{K_1^*,\ldots,K_{n_0}^*\}$, draw
$Z_k\vert\bmu,\btau$ from the prior multivariate normal density $
\phi(Z|\bmu,\btau)$. For $j=1,\ldots,n_0$, draw
$Z_{K_j^*}:={\mathbf
{\bbeta}}^*_j$ from the density proportional to
$\phi(\bbeta^*_j|\bmu,\btau)\prod_{\{i:K_i=K_j^*\}}L(Y_i,\bbeta^*_j)$
by using, for example, a standard Metropolis--Hastings procedure.
\item\textit{Conditional draw for $\theta=(\bmu,\btau)$.}
Conditional on
$\btau,\K,\Z,\Y$, draw $\bmu$ from a multivariate normal
distribution with parameters
\[
\frac{\lambda\m+n_0\bar{\bbeta}_{n_0}}{\lambda+n_0}
\quad\mbox{and}\quad\frac{\btau}{\lambda+n_0},
\]
where $\bar{\bbeta}_{n_0}=n_0^{-1}\sum_{j=1}^{n_0}\bbeta^*_j$.
Conditional on $\K,\Z,\Y$, draw $\btau$ from an inverse Wishart
distribution with parameters
\[
\nu_0+n_0\quad\mbox{and}\quad\frac{\nu_0\S_0+n_0\S_{n_0}+R({\bar{\bbeta}
}_{n_0},\m)}{\nu_0+n_0},
\]
where
\[
\S_{n_0}=\frac{1}{n_0}\sum_{j=1}^{n_0}
(\bbeta^*_j-{\bar{\bbeta}}_{n_0})(\bbeta^*_j-{\bar{\bbeta}}_{n_0}
)^\prime
\quad\mbox{and}\quad
R({\bar{\bbeta}}_{n_0} ,\m)=\frac{\lambda n_0}{\lambda+n_0}
({\bar{\bbeta}}_{n_0}-\m)({\bar{\bbeta}}_{n_0}-\m)^\prime.
\]
\end{enumerate}
\end{alg}

Notice that, when $n_{0}=1$, Steps 3 and 4 reduce to the MCMC steps
for a parametric Bayesian model. Iterating the steps above produces
a draw from the distribution $\Z,\K,\p,\theta\vert\Y$. Thus, each
iteration $m$ defines a probability measure $G^{(m)}(\cdot)=
\sum_{k=1}^{N}p_k^{(m)}\delta_{Z_k^{(m)}}(\cdot)$, which eventually\vspace*{-3pt}
approximates draws from the posterior distribution\vspace*{1pt} of $\tilde
G\vert\Y$. Consequently, one can approximate the posterior
distributional properties of the choice probabilities
$\rmP(\{j\}\vert\tilde G,\mathbf{x})$ by constructing (iteratively)
\[
\rmP\bigl(\{j\}\vert G^{(m)},\mathbf{x}\bigr)
=\sum_{k=1}^Np_k^{(m)} \frac{\exp\{\mathbf{x}_j^\prime Z_k^{(m)}\}}
{ \sum_{l\in\mathbf{C} }\exp\{\mathbf{x}_l^\prime Z_k^{(m)}\}};
\]
see \eqref{eq:BMMNL}. For instance, an histogram of the
$\rmP(\{j\}\vert G^{(m)},\mathbf{x})$, for $m=1,\ldots, M$,
approximates the
posterior distribution. An approximation to the posterior mean
$\mathrm{E}[\rmP(\{j\}\vert\tilde{G},\mathbf{x})\vert\Y]$
is obtained by $M^{-1}\sum_{m=1}^{M}\rmP(\{j\}\vert G^{(m)},\mathbf{x})$ or,
alternatively, by
%
\begin{equation}\label{eq:postmean}
\widehat{{P}} (\{j\}\vert\mathbf{x}):=\frac{1}{M}\sum_{m=1}^{M}\mathrm
{E} \bigl[\rmP(\{j\}\vert
\tilde{G},\mathbf{x})\vert\theta^{(m)},\bbeta_1^{(m)},\ldots,
\bbeta_n^{(m)} \bigr],
\end{equation}
where $\mathrm{E} [\rmP(\{j\}\vert
\tilde{G},\mathbf{x})\vert\theta,\bbeta_{1},\ldots, \bbeta_{n} ]$ is given
in \eqref{eq:prediction} and $\bbeta_i^{(m)}=Z^{(m)}_{K_{i}^{(m)}}$.

%
%
\section{Bayesian modeling for panel data}\label{section:4}
The MMNL framework may also be used to model choice probabilities
based on panel data. In the panel data setting, each individual $i$
is observed to make a sequence of choices at different time points.
The random utility for choosing ${j}$ for individual $i$ in choice
situation $t$ is given by
\[
U_{ijt}=\mathbf{x}^\prime_{ijt}\bbeta_i+\varepsilon_{ijt},\qquad j\in\mathbf{C},
\]
for times $t=1,\ldots, T_i$.
The MMNL model can be described as follows [see Train (\citeyear{train03}), Section~6.7]: given $\bbeta_{i}$, the probability that a person makes the
sequence of choices $\Y_i=\{Y_{i1},\ldots, Y_{iT_i}\}$
is the product of logit formulae
\[
L(\Y_i,\bbeta_i)= \prod_{t=1}^{T_i}\frac{\exp\{\mathbf{x}_{iY_{it}t}^\prime\bbeta_i\} }
{\sum_{j\in\mathbf{C} }\exp\{\mathbf{x}_{ijt}^\prime\bbeta_i\}}.
\]
The MMNL model is completed by taking the $\bbeta_{i}$ to be from
a distribution $G$ so that the unconditional choice probability is
specified by
\[
\mathrm{P}(\Y_{i}\vert G,\mathbf{x}_i)
=\int_{\rr^{d}} \prod_{t=1}^{T_i}\frac{\exp\{\mathbf{x}_{iY_{it}t}^\prime\bbeta\} }
{\sum_{j\in\mathbf{C} }\exp\{\mathbf{x}_{ijt}^\prime\bbeta\}}G(\ddr
\bbeta)
=\int_{\rr^{d}} L(\Y_i,\bbeta)G(\ddr\bbeta),\vadjust{\goodbreak}
\]
where $\mathbf{x}_i=\{\mathbf{x}_{ijt}, j\in\mathbf{C}, t=1,\ldots,T_i\}$
denotes the
array of covariates associated with the sequence of choices of
individual $i$. Similarly to the non-panel data setting, we wish to
model $G$ as a random probability measure in a Bayesian
framework. While it is possible to choose $\tilde{G}$ to follow a
Dirichlet process, this would result in possible ties among the
individual's preferences~$\bbeta_{i}$. In order to preserve the
distinct nature of each individual's preference, we assume that,
given $\tilde{G}$, the $\bbeta_i$ are i.i.d.~with distribution
$\tilde{G}$, where $\tilde{G}$ is a mixture of multivariate normal
distributions with random mixing distribution $\tilde\Pi$. That is,
$\tilde{G}$ has random density
$\int_{\Theta} \phi(\bbeta|\bmu,\btau)
\tilde\Pi(\ddr\bmu,\ddr\btau)$, where
$\Theta=\rr^d\times\mathcal{S}$ with $\mathcal{S}$ the space of
covariance matrices.
Specifically, we take $\tilde\Pi$ to be a Dirichlet process with
shape $aF$, $F$ a probability measure on $\Theta$. Hence, the
Bayesian MMNL model for individual $i$ is expressible as
\[
\mathrm{P}(\Y_{i}\vert \tilde{G},\mathbf{x}_i)=
\int_{\rr^d}L(\Y_i,\bbeta)\tilde{G}(\ddr\bbeta)
=\int_{\rr^d}\int_{\Theta}L(\Y_i,\bbeta)\phi(\bbeta|\bmu,\btau)
\tilde\Pi(\ddr\bmu,\ddr\btau)\,\ddr\bbeta.
\]
While one may use any choice for $F$, we take
$F(\ddr\bmu,\ddr\btau)$ to be the multivariate
normal inverse Wishart distribution $\operatorname{N\mbox{-}IW}(\m,{\lambda}^{-1}\btau,\S_{0},\nu_{0})$ described in
Section~\ref{section:3}.


\subsection{Blocked Gibbs algorithm for panel data}
The explicit posterior analysis for the panel data case is quite
similar to the non-panel case. The main difference is that the
$(\bmu_i,\btau_i)$, $i=1,\ldots, n$, rather than $\bbeta_1,\ldots,
\bbeta_n$, are drawn from the Dirichlet process. Here, we will
briefly focus on the relevant data structure and then proceed to a
description of how to implement the blocked Gibbs sampler.
The joint distribution of the augmented data can be expressed using
a hierarchical model as follows:
\begin{eqnarray}\label{eq:panelmodel2}
\Y_i\vert \bbeta_{i}&\ind&L(\Y_i,\bbeta_i)
\qquad\mbox{for }i=1,\ldots,n\mbox{ and }Y_{it}\in\mathbf{C},\nonumber\\
\bbeta_i\vert \bmu_i,\btau_i &\ind&\phi(\bbeta_{i}|\bmu
_i,\btau_i)\qquad
\mbox{for }i=1,\ldots,n,\nonumber\\[-8pt]\\[-8pt]
\bmu_i,\btau_i\vert\tilde\Pi&\iid&\tilde\Pi\qquad
\mbox{for }i=1,\ldots,n, \nonumber\\
\tilde\Pi&{\sim}&\pp(\ddr\Pi|a F). \nonumber
\end{eqnarray}
Similar to the non-panel case, the blocked Gibbs sampler works by
using the $\pp^{\ssN}(\ddr\Pi|a F)$ in place of the law of the
Dirichlet process $\pp(\ddr\Pi|a F)$. We now sample
$(\K,\p,\Z,\bbeta_1,\ldots,\bbeta_n)$ from the distribution
proportional to
\[
\Biggl[\prod_{i=1}^nL(\Y_i,\bbeta_i)\phi(\bbeta_i|\bmu_i,\btau_i) \Biggr]
\Biggl[\prod_{i=1}^n\sum_{k=1}^Np_k\delta_{Z_k}(\ddr\bmu_i,\ddr\btau
_i) \Biggr]
\pi(\p)\prod_{k=1}^{N}F(\ddr Z_{k}).
\]
Here, we use the fact that $(\bmu_{i},\btau_{i})=Z_{K_{i}}$ for
$i=1,\ldots, n$. To approximate the posterior law of various
functionals, we cycle through the following steps.\vadjust{\goodbreak}

\begin{alg}\label{alg2}
\begin{enumerate}
\item\textit{Conditional draw for $\K$.}
Independently sample $K_i$ according to
\[
\rmP\{K_i\in\cdot\vert\p,\Z,\bbeta_1,\ldots,\bbeta_n,\Y\}
=\sum_{k=1}^N
p_{k,i}\delta_k(\cdot)\qquad\mbox{for } i=1,\ldots,n,
\]
where $(p_{1,i},\ldots,p_{N,i})\propto (p_1
\phi(\bbeta_{i}|Z_1),\ldots,p_N \phi(\bbeta_i|Z_N) )$.
\item\textit{Conditional draw for $\p$.} $p_1=V_1^*$,
$p_k=(1-V_1^*)\cdots(1-V_{k-1}^*)V_k^*$, $k=2,\ldots,N-1$ and
$V_N^*=1$, where, if $e_k$ records the number of $K_i$ values which
equal $k$,
\[
V_k^*\ind\operatorname{beta} \Biggl(1+e_k,a+\sum_{l=k+1}^{N}e_l
\Biggr),\qquad k=1,\ldots,N-1.
\]
\item\textit{Conditional draw for $\Z$.} Let
$\{K_1^*,\ldots,K_{n_0}^*\}$ denote the unique set of $K_i$ values.
For each $k\notin\{K_1^*,\ldots,K_{n_0}^*\}$, draw
$Z_k=(\bmu_{k},\btau_{k})$ from the prior
$\operatorname{N\mbox{-}IW}(\m,{\lambda}^{-1}\btau,\S_{0},\nu_{0})$. For
$j=1,\ldots,n_0$, draw $Z_{K_j^*}:=(\bmu^{*}_j,\btau^{*}_j)$ as
follows: (a) conditional on $\btau^{*}_{j},\K,\bbeta_{1}, \ldots,\break
\bbeta_{n},\Y$, draw $\bmu^{*}_j$ from a multivariate normal
distribution with parameters
\[
\frac{\lambda\m+e_{K_j^*}\bar{\bbeta}^*_j}{\lambda
+e_{K_j^*}}\quad\mbox{and}\quad\frac{\btau^*_j}{\lambda+e_{K_j^*}},
\]
where
$\bar{\bbeta}^*_j=(e_{K_j^*})^{-1}\sum_{\{i:K_i=K_j^*\}}\bbeta_i$;
(b) conditional on $\K,\bbeta_{1},\ldots,\bbeta_{n},\Y$, draw
$\btau^*_j$ from an inverse Wishart distribution with parameters
\[
\nu_{0}+e_{K_j^*}\quad\mbox{and}\quad\frac{\nu_0\S_0+e_{K_j^*}\S_j
+R(\bar{\bbeta}^*_j,\m)}{\nu_0+e_{K_j^*}},
\]
where
\[
\S_j=\frac{1}{e_{K_j^*}}\sum_{\{i:K_i=K_j^*\}}
(\bbeta_i-\bar{\bbeta}^*_j)(\bbeta_i-\bar{\bbeta}^*_j)^\prime
\quad\mbox{and}\quad
R(\bar{\bbeta}^*_j,\m)=\frac{\lambda e_{K_j^*}}{\lambda+e_{K_j^*}}
(\bar{\bbeta}^*_j-\m)(\bar{\bbeta}^*_j-\m)^\prime.
\]
\item
\textit{Conditional draw for $\bbeta_1,\ldots, \bbeta_n$.} For each
$j=1,\ldots, n_{0}$, independently draw $\bbeta_i$,
$i\in\{l\dvt K_{l}=K_j^{*}\}$, from the density proportional to
$L(\Y_i,\bbeta_i)\phi(\bbeta_{i}|\bmu^{*}_{j},\btau^{*}_{j})$ by
using, for example, a standard Metropolis--Hastings procedure.
\end{enumerate}
\end{alg}

When $n_{0}=1$, Steps 3 and 4 equate with a parametric MCMC procedure
for panel data models similar to the algorithm described in
Train (\citeyear{train03}), Section 12.


\section{Simulation study}\label{section:5}
In this section, we present some empirical evidence that shows how
the MMNL procedures perform overall and relative to GML models and
finite mixture (FM) of MNL models. We proceed to the estimation of the
choice probabilities
based on simulated data. Two different artificial data sets are
generated for the simulation study: data set 1 is produced for studying
non-panel data models, while data set 2 is designed to study models
with panel data.
In both cases, we consider a RUM model with three possible responses
($J=3$) relative to the utilities $U_{1},U_{2}$ and $U_{3}$,
\[
\cases{
U_1=x_{11}\beta_1+x_{12}\beta_2+\varepsilon_1,\cr
U_2=x_{21}\beta_1+x_{22}\beta_2+\varepsilon_2,\cr
U_3=x_{31}\beta_1+x_{32}\beta_2+\varepsilon_3.
}
\]
As for data set 1, we choose $\varepsilon_1,\varepsilon_2,
\varepsilon_3\iid \mathrm{standard}$ Gumbel and
$\bbeta=(\beta_1,\beta_2)^\prime\iid 0.5 \times\delta_{(-5,5)}
+ 0.5 \times\delta_{(5,-5)}$.
For individual $i$, we randomly generate (componentwise) the
covariate vector ${\mathbf x}_{i}=(x_{11}, x_{12}, x_{21}, x_{22},
x_{31},x_{32})$, independently from a $\operatorname{Uniform}( -2,2)$ distribution.
Set $Y_{i}=j$ if $U_{ij}>U_{il}$, $l\neq j$, for $j=1,2,3$. Repeat
this procedure $n$ times independently to obtain a data set with
$(Y_i,\mathbf{x}_i)$ for $i=1,\ldots,n$. As for data set 2, we
assume that there are $n$ individuals, each making $T_{i}=10$ choices for
$i=1,\ldots,n$.
We then simulate data using the same model used to generate data set
1. The only change is that $\bbeta$ is drawn from the two-component
mixture of bivariate normal distributions,
$\bbeta\iid 0.5\times N((-5,5)^{\prime},2\mathbf{I})+0.5\times
N((5,-5)^{\prime},2\mathbf{I})$, where $\mathbf{I}$ is the identity
matrix.

We start by applying our procedures to the estimation of choice
probabilities $\rmP(\{j\}\vert \mathbf{x})$, for $j=1, 2, 3$, based on the
set of covariates $\mathbf{x}=(1.0, -0.9, 1.0, 0.2, 1.0, 0.9)$.
The prior parameters for the specifications of the Bayesian MMNL
models for panel and non-panel data (pertaining to the explicit
models in Sections \ref{section:3} and \ref{section:4}) are set to
be $a=1$, $\nu_0=2$, $\m=( 0,0) ^\prime$, $\S_0=\mathbf{I}$ and
$\lambda=1$. Additionally, we use $N=100$ for the truncation level in
the blocked Gibbs Algorithms \ref{alg1} and \ref{alg2} given in Sections \ref{section:3}
and \ref{section:4}, respectively.
A Bayesian GML model is also estimated for comparison with the
same specifications for $\nu_{0}$, $\mathbf{m}$, $\S_{0}$ and
$\lambda$. In all cases, we use the estimator \eqref{eq:postmean}
based on an initial burn-in of 10,000 cycles and an additional 10,000
Gibbs cycles ($M=10\mbox{,}000$)
for the estimation.
In addition, to measure how good of our estimates are, we define a
measure, root mean square (RMS) value, as
\[
\mathrm{RMS} = \sqrt{\frac{1}{J}\sum_{j\in\mathbf{C}}
\frac{1}{ M}\sum_{m=1}^M \bigl(\rmP\bigl(\{j\}\vert G^{(m)},\mathbf{x}\bigr)-
\rmP_0(\{j\}\vert\mathbf{x}) \bigr)^2},
\]
where $\rmP_0(\{j\}\vert\mathbf{x})$ is the choice probability resulting
from the data generating process.

\begin{table}
\tabcolsep=0pt
\caption{Simulation results for data set 1 (columns 3--4) and for data
set 2 (columns 5--6) with $\x=(1.0,-0.9,1.0,0.2,1.0,0.9) $ -- the
estimates (Est.), the credible intervals (C.I.) and the root mean
square (RMS) values are presented; GML${}={}$Gaussian mixed logit, MMNL${}={}$mixed multinomial logit}\label{table:1}
\begin{tabular*}{\textwidth}{@{\extracolsep{\fill}}llllllllll@{}}
\hline
&   &\multicolumn{4}{l}{Data set 1 (non-panel case)}&\multicolumn{4}{l@{}}{Data set 2 (panel case)}\\
&   &\multicolumn{4}{l}{$n=500$}&\multicolumn{4}{l@{}}{$n=100$, $T_i=10$}\\[-6pt]
&&\multicolumn{4}{c}{\hrulefill}&\multicolumn{4}{c@{}}{\hrulefill}\\
&                       & True  & Est. & (95\% C.I.) & RMS & True  & Est.& (95\% C.I.) & RMS    \\
\hline
GML &$\rmP(\{1\}\vert  \x)$ & 0.4980 & 0.3203&(0.2907, 0.3501) &  & 0.4939 & 0.4585&(0.4476, 0.4685) &  \\
&$\rmP(\{2\}\vert  \x)$ & 0.0167 & 0.3348&(0.3308, 0.3377) &  & 0.0279 & 0.0521&(0.0378, 0.0675) &  \\
&$\rmP(\{3\}\vert  \x)$ & 0.4853 & 0.3449&(0.3191, 0.3715) &  & 0.4782 & 0.4894&(0.4717, 0.5061) &  \\
&                    &       &      & & 0.2258 &&      &  &  0.0266      \\[6pt]
MMNL&$\rmP(\{1\}\vert  \x)$ & 0.4980 & 0.4856&(0.4748, 0.4945) &  & 0.4939 & 0.4586&(0.4495, 0.4670) &  \\
&$\rmP(\{2\}\vert  \x)$ & 0.0167 & 0.0257&(0.0069, 0.0551) &  & 0.0279 & 0.0494&(0.0329, 0.0679) &  \\
&$\rmP(\{3\}\vert  \x)$ & 0.4853 & 0.4886&(0.4615, 0.5057) &  & 0.4782 & 0.4920&(0.4705, 0.5107) &  \\
&                    &       &       &      & 0.0137  &&   &    &0.0265             \\
\hline
\end{tabular*}
\end{table}

\begin{figure}

\includegraphics{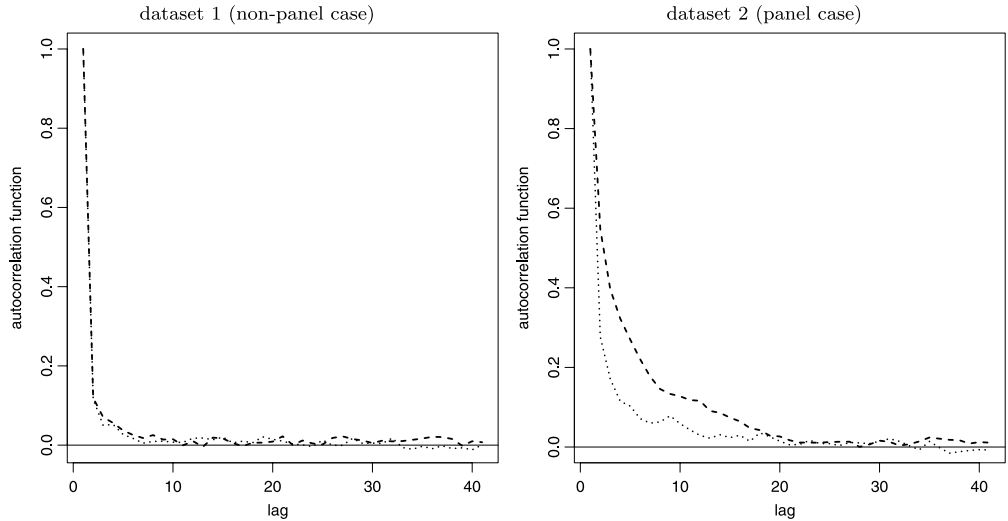}

  \caption{MMNL model: Autocorrelation functions for the choice probability
  $\mathrm{P}(\{1\}\vert  \x)$ for data set 1 (left) and data set 2
(right), obtained from the posterior sample of the $\bbeta$'s for
the MMNL model with prior hyperparameter $\lambda=0.01$ (dashed)
and $\lambda=1$ (dotted).}\label{figure:1}
\end{figure}

Simulation results using data set 1 ($n=500$) and data set 2 ($n=100,
T_i=10$) are summarized in Table \ref{table:1}, together with RMS
values, for both the GML and the MMNL models. They show that the
performance of the nonparametric MMNL model is better than that of the
parametric GML model in the non-panel case, as indicated by a smaller
RMS value and more accurate estimates of choice probabilities, while
the GML and MMNL models display similar performances in the panel case.
As expected, the GML model suffers from misspecification in the
non-panel case, while the two-component mixture of bivariate normals
used for generating data set 2 is correctly accounted for by the GML
because of the hyperprior on the parameter $(\bmu,\btau)$ we are
using. We then get confirmation that the fit of the MMNL model is as
good as that of the GML model.
We also performed estimation of the MMNL model for different choices of
the scale parameter $\lambda$ (not reported here) which show two
different behaviors for the non-panel and the panel case. As for the
non-panel case, RMS values and the estimates remain stable, whereas, in
the panel case, the estimates are more accurate when we decrease
$\lambda$ with slightly smaller RMS values. An interpretation of an
increase of accuracy is as follows: a smaller $\lambda$ corresponds to
a more diffuse $H$, the
prior predictive distribution of $\tilde G$. Since $H$ is different
from the distribution used to simulate the $\bbeta$'s in the data
generating process, we obtain evidence that a diffuse $H$ helps in
capturing the true form of the mixing distribution $G$. Also, note
that a smaller $\lambda$ yields a smaller RMS, the latter being a
measure of the combination of the accuracy and the variability of
the posterior variates of $\rmP(\{j\}\vert \mathbf{x})$.
An examination of their autocorrelation functions along the chain
shows that a smaller $\lambda$
causes a slower mixing of the blocked Gibbs sampler, which increases the
component of variability in the RMS; see Figure \ref{figure:1}. The
decrease in RMS then shows that such precision loss is more than
balanced by a higher accuracy of the estimate, although one should
also control the convergence properties of the sampler by
avoiding taking $\lambda$ too small.

We investigated the sensitivity of the results to the prior
parameter $\nu_0$, where a larger $\nu_0$ corresponds to a more
concentrated inverse Wishart distribution on $\S_0$. However, we did not
observe substantial differences in the estimation by varying $\nu_0$
and we decided to set $\nu_0=2$ and $\S_0=\mathbf{I}$ as a default
non-informative choice for these parameters; see Train (\citeyear{train03}), Section
12. The nonparametric prior on $\tilde{G}$ is also dependent on the total
mass $a$, which is positively related to the number of components in
the mixture distribution of the $\bbeta$'s. Generally, $a=1$ is
considered a default choice for a finite mixture model with a fixed,
but uncertain, number of components. We performed estimation for
larger $a$, obtaining almost identical results: $a=1$ was, in fact,
sufficient for detecting the two-component mixture we used in
generating the data. Although we have not done so, the blocked Gibbs
procedures described in Sections \ref{section:3} and \ref{section:4}
can be easily extended to place an additional prior on $a$.
Furthermore, the truncation level of $N=100$ in \eqref{eq:dpn} is
sufficiently large as we observed almost identical estimation
results from runs of the blocked Gibbs sampler with larger values of $N$.

\begin{table}[b]
\caption{MMNL model: estimates and the root mean square (RMS) for
data set 1 and for data set 2 with $\x=(1.0,-0.9,1.0,0.2,1.0,0.9)$ and
different sample sizes}\label{table:2}
\begin{tabular*}{\textwidth}{@{\extracolsep{\fill}}lllllllll@{}}
\hline
&\multicolumn{4}{l}{Data set 1 (non-panel case)}&
\multicolumn{4}{l@{}}{Data set 2 (panel case)}\\[-5pt]
&\multicolumn{4}{c}{\hrulefill}&\multicolumn{4}{c@{}}{\hrulefill}\\
&&&&&&$n=10$ & $n=50$ & $n=100$ \\
& True & $n=50$ & $n=100$ & $n=500$ & True &  $T_i=10$ & $T_i=10$&$T_i=10$\\
\hline
$\rmP(\{1\}\vert \x)$ & 0.4980 & 0.4927 & 0.5145 & 0.4856 &  0.4939 & 0.5956 & 0.4176 & 0.4586 \\
$\rmP(\{2\}\vert \x)$ & 0.0167 & 0.1046 & 0.0489 & 0.0257 &  0.0279 & 0.0527 & 0.0562 & 0.0494 \\
$\rmP(\{3\}\vert \x)$ & 0.4853 & 0.4027 & 0.4366 & 0.4886 &  0.4782 & 0.3517 & 0.5261 & 0.4920 \\[3pt]
RMS &    & 0.0867  & 0.0440  & 0.0137 && 0.0977  & 0.0556 & 0.0265    \\
\hline
\end{tabular*}
\end{table}

\begin{figure}

\includegraphics{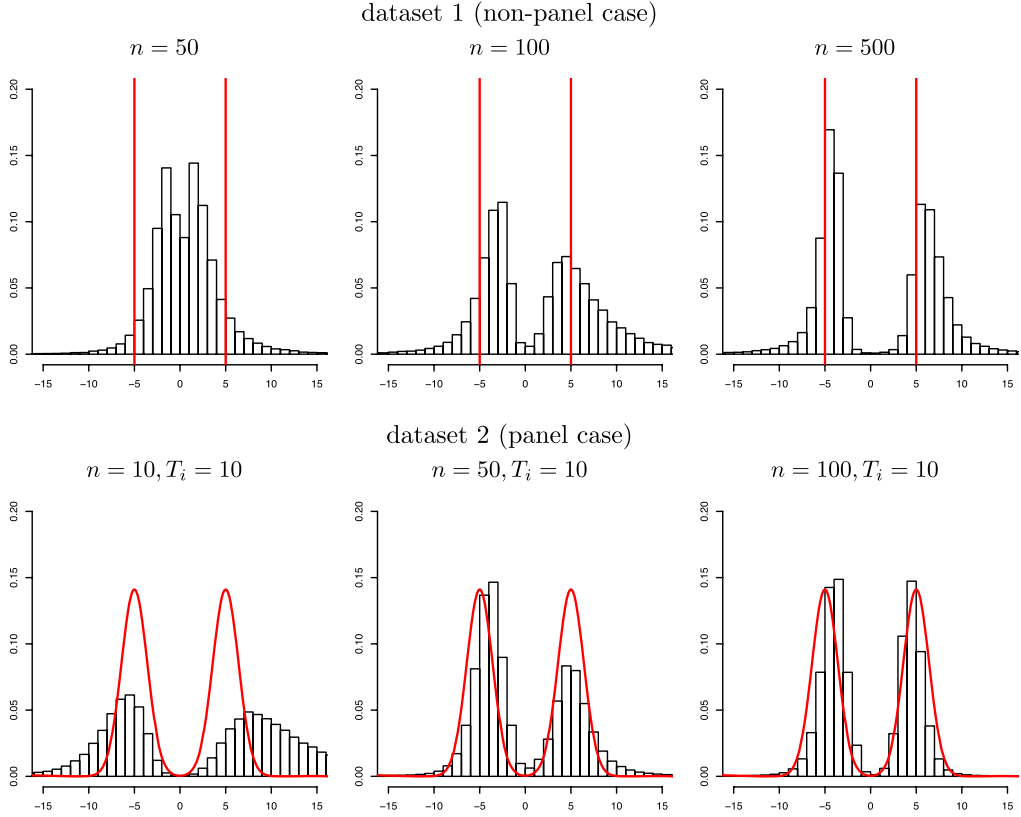}

\caption{MMNL model: histogram estimate of the posterior marginal density of $\beta_1$'s for
data set 1 (top) and for data set 2
(bottom) and different sample sizes. The solid lines represent the true mixing distribution.}\label{figure:2}
\vspace*{-6pt}
\end{figure}

The second simulation study aims at the verification of the consistency
result of Section \ref{section:2} by estimating the MMNL model for increasing
sample sizes for both data set 1 and data set 2. We also sample $\bbeta
$ variates from their posterior distribution, thus obtaining
approximated evaluation of the mixing distribution $G$.
The prior parameters are set as $a=1$, $\nu_0=2$, $\m=( 0,0)
^\prime$, $\S_0=\mathbf{I}$, $N=100$ and $\lambda=1$. Table
\ref{table:2} reports the results by showing, as expected, a noticeable
decrease of RMS for both non-panel and panel data as the number of
observations increases. In addition, Figure \ref{figure:2} reports
the histograms of samples for $\beta_1$ from its marginal posterior
distribution against the mixing distribution used in the data
generating process: it shows how the approximation of the true
mixing distribution $G$ improves as more and more data become
available.

Finally, we evaluate the performance of the Bayesian MMNL model via a
comparison with the finite mixture (FM) MNL model estimated via the EM
algorithm described in Train (\citeyear{train08}), Section 4. The FM MNL model can be
considered nonparametric in the sense that the locations and weights of
the mixing distribution $G$ are both assumed to be parameters. The
selection of the number of points in the mixing is based on the BIC
information criterion. We consider 500 Monte Carlo replicates of each
of the following 6 situations: data set 1 with sample sizes $n=50,100$
and $500$; data set 2 with $(n=10,T_i=10)$, $(n=50,T_i=10)$ and
$(n=100,T_i=10)$.
For a given sample, the posterior estimate of $P(\{j\}\vert \mathbf{x})$
in equation (\ref{eq:postmean}) is computed, based on $6000$ Gibbs cycles after a
burn-in period of 4000 for $j=1,2,3$ and for $\mathbf{x}$ in a
$6$-dimensional grid of the hypercube $(-2,2)^6$ of $5^6$
equally-spaced points. At the same time, we compute the FM MNL estimate
of $P(\{j\}\vert \mathbf{x})$ for $j=1,2,3$, evaluated on the same grid
of $\mathbf{x}$-points. We call $\hat\q(\mathbf{x})$ and $\q_0(\mathbf{x})$
the estimated vector and the true vector of choice probabilities
evaluated at $\mathbf{x}$, respectively.
%
\begin{table}[h!]
\caption{Average $L_1$-error from 500 Monte Carlo replicates -- FM MNL $=$ finite mixture of
multinomial logit; MMNL $=$ mixed multinomial logit}\label{table:3}
\begin{tabular*}{\textwidth}{@{\extracolsep{\fill}}lllllll@{}}
\hline
&\multicolumn{3}{l}{Data set 1 (non-panel case)}&\multicolumn{3}{l@{}}{Data set 2 (panel case)}\\[-5pt]
&\multicolumn{3}{c}{\hrulefill}&\multicolumn{3}{c@{}}{\hrulefill}\\
&&&& $n=10$  & $n=50$  & $n=100$ \\
& $n=50$ & $n=100$ & $n=500$ & $T_i=10$& $T_i=10$ & $T_i=10$\\
\hline
FM MNL & 0.0521 & 0.0295 & 0.0107 & 0.0891 & 0.0505 & 0.0297 \\
MMNL   & 0.0577 & 0.0316 & 0.011\phantom{0}  & 0.0827 & 0.0467 & 0.0268 \\
\hline
\end{tabular*}
\vspace*{-6pt}
\end{table}
We measure the overall error of
estimation with the $L_1$-distance $\int_{\mathcal{X}}|\hat\q(\mathbf{x})-\q
_0(\mathbf{x})|\,\ddr\mathbf{x}$, which corresponds to the (rescaled) distance
$d(\hat\q,\q_0)$ in equation (\ref{eq:hellinger}), with $M(\ddr\mathbf{x})$ being the
uniform distribution on the hypercube $(-2,2)^6$. We compute the
$L_1$-error for the Bayesian MMNL estimator and the FM MNL estimator,
then average over the 500 Monte Carlo replicates.
 The results are
reported in Table \ref{table:3} and show that the MMNL estimators outperform the FM
MNL estimators in the panel case for all sample sizes, while in the
non-panel case, the situation is reversed, with a similar performance
for $n=500$.
Note, however, that data set 1 is generated exactly from a finite
mixture model so that the FM MNL model is expected to perform well.
Overall, the decrease in the average error for larger sample sizes is a
further confirmation of the consistency result of Section 2.

\section{\texorpdfstring{Proof of Theorem \protect\ref{theorem:consistency}}{Proof of Theorem 1}}
\label{section:6}
Throughout this section, we work with the family of multinomial
logistic kernels
\[
k_j(\mathbf{x},\bbeta)=\frac{\exp(x_{j}^\prime\bbeta)}{
\sum_{l\in\mathbf{C}}\exp(x_{l}^\prime\bbeta)},\qquad j=1,\ldots,J.
\]
With $q_j(\mathbf{x};G)$ denoting the $j$th element of the vector
$\q(\mathbf{x};G)$, we have that
$q_j(\mathbf{x};G)=\int_{\rr^d}k_j(\mathbf{x},\bbeta)G(\ddr\bbeta)$.
Note that $q_{_Y}(\mathbf{x};G_0)$ is the joint density of $(Y,\X)$ with
respect to the counting measure on the integer set $\mathbf{C}$ and the
measure $M(\ddr\mathbf{x})$ on $\mathcal X$.

For the proof of Theorem \ref{theorem:consistency}, the following
lemma is essential, stating that on the space $\mathbb{P}$, the
weak topology and the topology induced by the $\ll_1$-distance $d$
defined in \eqref{eq:hellinger} are equivalent.
%
%
\begin{lemma}\label{lemma:1}
Let $d_w$ be any distance that metrizes the weak topology on
$\mathbb{P}$ and $(G_n)_{n\geq1}$ be a sequence in $\mathbb{P}$.
Then
$d_w(G_n,G_0)\to0$ if and only if
$d (\q(\cdot;G_n),\q(\cdot;G_0) )\to0$.
\end{lemma}
\begin{pf}
For the ``if'' part, it is sufficient that $d_w(G_n,G_0)\to0$ implies that
$\int_{\mathcal X} |q_j(\mathbf{x};G_n)-
q_j(\mathbf{x};G_0) |M(\ddr\mathbf{x})\to0$
for an arbitrary $j\in\mathbf{C}$. The latter is a consequence of the
definition of weak convergence and an application of Scheff\'{e}'s
theorem since $k_j(\mathbf{x},\bbeta)$ is bounded and continuous in
$\bbeta$ for each $\mathbf{x}\in\mathcal X$.
To show the converse, we prove that $G $ being distant from $G_0$ in the
weak topology implies that $\q(\cdot;G)$ is distant from
$\q(\cdot;G_0)$ in the $\ll_1$-distance $d$. Define a weak
neighborhood of $G_0$ as
\[
V= \biggl\{G\dvt
\biggl|\int_{\rr^d}\int_{\mathcal X}k_j(\mathbf{x},\bbeta)M(\ddr\mathbf{x})G(\ddr\bbeta)-
\int_{\rr^d}\int_{\mathcal X}k_j(\mathbf{x},\bbeta)M(\ddr\mathbf{x})G_0(\ddr\bbeta) \biggr|<\delta,
j\in\mathbf{C} \biggr\}.
\]
Since $\int_{\mathcal X}k_j(\mathbf{x},\bbeta)M(\ddr\mathbf{x})$ is a bounded
continuous function on $\rr^d$ for each $j$, $G\in V^c$ implies that
$d_w(G,G_0)>\delta$. Based on the inequalities
\begin{eqnarray*}
d (\q(\cdot;G),\q(\cdot;G_0) )
&\geq&\max_{j\in\mathbf{C}}
\int_{\mathcal X} |q_j(\mathbf{x};G_n)-q_j(\mathbf{x};G_0) |M(\ddr\mathbf{x})\\
&\geq&\max_{j\in\mathbf{C}}
\biggl|\int_{\mathcal X}\int_{\rr^d}k_j(\mathbf{x},\bbeta)G(\ddr\bbeta
)M(\ddr\mathbf{x})-
\int_{\mathcal X}\int_{\rr^d}k_j(\mathbf{x},\bbeta)G_0(\ddr\bbeta
)M(\ddr\mathbf{x}) \biggr|
\end{eqnarray*}
and an application of Fubini's theorem, it follows that, for any
$\epsilon<\delta$ and any $G\in V^c$,
$d (\q(\cdot;G),\q(\cdot;G_0) )>\epsilon$. The proof is
then complete.
\end{pf}

\begin{remark}\label{remark:1}
Lemma \ref{lemma:1} has two important consequences: (a) both
$\mathscr{Q}$ and $\mathbb{P}$ are separable spaces under the metric
$d$; (b) the statement of Theorem \ref{theorem:consistency} is
equivalent to saying that $\pp_n$ accumulates all probability mass
in a weak neighborhood of $G_0$.
\end{remark}

Define $\Lambda_n(G)=\prod_{i=1}^n q_{_{Y_i}}(\X_i;G)
/q_{_{Y_i}}(\X_i;G_0)$ so that the posterior distribution of
$\tilde{G}$ can be written as
%
\begin{equation}\label{eq:Pi(A)}
\pp_n(A)=
\frac{\int_A\Lambda_n(G)\pp(\ddr G)}{\int_{\mathbb{P}} \Lambda
_n(G)\pp(\ddr G)}.
\end{equation}
We now take $A=\{G\dvt
d (\q(\cdot;G),\q(\cdot;G_0) )>\epsilon\}$ and will,
as is usual in the Bayesian consistency literature,
separately consider the numerator and the denominator of \eqref
{eq:Pi(A)}. To
this end, define $I_n=\int_{\mathbb{P}}\Lambda_n(G)\pp(\ddr G)$.
Relying on the separability of $\mathbb{P}$ under the topology
induced by $d$ (see Remark \ref{remark:1}), for any $\eta>0$, we can
cover $A$ with a countable union of disjoint sets $A_j$ such that
%
\begin{equation}\label{eq:cover}
A_j\subseteq A_j^*=\{G\dvt d (\q(\cdot;G),\q(\cdot;G_j) )<\eta\}
\end{equation}
and $\{G_j\}_{j\geq1}$ is a countable set in $\mathbb{P}$ such that
$d (\q(\cdot;G_j),\q(\cdot;G_0) )>\epsilon$ for any $j$.
Consider the fact that
\[
\pp_n(A)=\sum_{j\geq1}\pp_n(A_j)
\leq\sum_{j\geq1}\sqrt{\pp_n(A_j)}=\sum_{j\geq1}\sqrt
{I_n^{-1}\int_{A_j}
\Lambda_n(G)\pp(\ddr G)}.
\]
Hence, Theorem \ref{theorem:consistency} holds if we prove that, for
all large $n$,
%
\begin{eqnarray}\label{eq:den}
\forall c>0,\qquad I_n&>&\exp(-nc)\qquad\mbox{a.s.}
\\
\label{eq:num}
\exists b>0\mbox{:}\qquad \sum_{j\geq1}\sqrt{\int_{A_j}
\Lambda_n(G)\pp(\ddr G)}&<&\exp(-nb)\qquad\mbox{a.s.}
\end{eqnarray}
As for \eqref{eq:den}, consider the Kullback--Leibler (KL) support
condition of $\pp$ defined by
%
\begin{equation}\label{eq:KL}
\pp\biggl\{G\dvt \int_{\mathcal{X}} K(G_0,G|\mathbf{x})M(\ddr\mathbf{x})<\epsilon
\biggr\}>0\qquad \forall\epsilon>0,
\end{equation}
where $K(G_0,G|\mathbf{x})=\sum_{j\in\mathbf{C}}q_j(\mathbf{x};G_0)\log
[{q_j(\mathbf{x};G_0)
/q_j(\mathbf{x};G)}]$. If $\pp$\vspace*{1pt} satisfies condition \eqref{eq:KL}, then
\eqref{eq:den} holds. To see this, it is sufficient to note that the
KL divergence of $q_{_Y}(\X;G)$ from $q_{_Y}(\X;G_0)$ with respect
to the measure $M(\ddr\mathbf{x})$ on $\mathcal{X}$ and the counting measure on
$\mathbf{C}$ is given by $\int K(G,G_0|\mathbf{x})M(\ddr\mathbf{x})$.
By the compactness of $\mathcal X$, the law of large numbers then
leads to
\[
\frac{1}{ n}\sum_{i=1}^n\log\frac{q_{_{Y_i}}(\X_i;G_0)
}{ q_{_{Y_i}}(\X_i;G)}\to
\int_{\mathcal X} K(G_0,G\vert\mathbf{x})M(\ddr\mathbf{x})\qquad
\mbox{a.s.}
\]
The result in \eqref{eq:den} then follows from standard arguments, see,
for example, Wasserman (\citeyear{wass98}).
Lemma \ref{lemma:2} below states that \eqref{eq:KL} is satisfied
under the hypotheses of Theorem \ref{theorem:consistency}.
%
%
\begin{lemma}\label{lemma:2}
If $G_0$ lies in the weak support of $\pp$ and condition \textup{(i)} of
Theorem \ref{theorem:consistency} holds, then $G_0$ is in the KL
support of $\pp$, according to \eqref{eq:KL}.
\end{lemma}
\begin{pf}
It is sufficient to show that for any $j\in\mathbf{C}$ and any
$\eta<1$, there exists a $\delta$ such that $|q_j(\mathbf{x};G)/
q_j(\mathbf{x};G_0)-1|\leq\eta$ whenever $G$ is in $W_\delta$, a
$\delta$-weak neighborhood of $G_0$. In fact, this implies that
\begin{eqnarray*}
\int_{\mathcal X}q_j(\mathbf{x};G_0)\log\biggl[\frac{q_j(\mathbf{x};G_0)
}{ q_j(\mathbf{x};G)} \biggr]M(\ddr\mathbf{x})
&\leq&\int_{\mathcal X}q_j(\mathbf{x};G_0) \biggl|\frac{q_j(\mathbf{x};G_0)
}{ q_j(\mathbf{x};G)}-1 \biggr|M(\ddr\mathbf{x})\\
&\leq&\int_{\mathcal X}q_j(\mathbf{x};G_0) \biggl(\frac{\eta}{1-\eta} \biggr)M(\ddr
\mathbf{x})\\
&\leq&\frac{\eta}{1-\eta},
\end{eqnarray*}
which, in turn, leads to the thesis, by the arbitrary nature of $j$.

Let $c=\inf_{\mathbf{x}\in\mathcal X}q_j(\mathbf{x};G_0)$, which is
positive by
condition (i) of Theorem \ref{theorem:consistency}, and assume
that $G\in W_\delta$ for a $\delta$ that will be determined later. Note
that, for any $\rho>0$, one can set $M_\rho>0$ such that
$G_0\{\bbeta\dvt |\bbeta|>M_\rho-\delta\}<\rho$. Then, using the
Prokhorov metric, $G\in W_\delta$ implies that $G\{\bbeta\dvt
|\bbeta|>M_\rho\}<\rho+\delta$. Also, note that
the family of functions $\{k_j(\mathbf{x},\bbeta), \mathbf{x}\in\mathcal X\}
$, as
$\bbeta$ varies in the compact set $\{|\bbeta|\leq M_\rho\}$, is
uniformly equicontinuous. By an application of the Arzel\`{a}--Ascoli
theorem, we know that, given a $\gamma>0$, there exist finitely many
points $\mathbf{x}_1,\ldots,\mathbf{x}_m$ such that, for any $\mathbf{x}\in
\mathcal X$,
there is an index $i$ such that
%
\begin{equation}\label{eq:Arz-Asc}
\sup_{|\bbeta|\leq M_\rho} |k_j(\mathbf{x},\bbeta)-
k_j(\mathbf{x}_i,\bbeta) |<\gamma.\vadjust{\goodbreak}
\end{equation}
For an arbitrary $\mathbf{x}\in\mathcal X$, choose the appropriate $\mathbf{x}_i$
such that \eqref{eq:Arz-Asc} holds, so that
\begin{eqnarray*}
\biggl|\frac{q_j(\mathbf{x};G)}{ q_j(\mathbf{x};G_0)}-1 \biggr|&\leq&\frac{1}{c} \biggl( \biggl|\int k_j(\mathbf{x}_i,\bbeta)G(\ddr\bbeta)
-\int k_j(\mathbf{x}_i,\bbeta)G_0(\ddr\bbeta) \biggr|\\
&&\hphantom{\frac{1}{c} \biggl(}
{}+
\int|k_j(\mathbf{x},\bbeta)-k_j(\mathbf{x}_i,\bbeta) |G(\ddr\bbeta)+
\int|k_j(\mathbf{x},\bbeta)-k_j(\mathbf{x}_i,\bbeta) |G_0(\ddr\bbeta) \biggr)\\
&:=&
\frac{I_1+I_2+I_3}{ c}.
\end{eqnarray*}
We have that $G\in W_\delta$ implies $I_1\leq\delta$.
As for $I_2$, we have
\begin{eqnarray*}
I_2&=&\int_{|\bbeta|\leq M_\rho} |k_j(\mathbf{x},\bbeta)-
k_j(\mathbf{x}_i,\bbeta) |G(\ddr\bbeta)
+\int_{|\bbeta|> M_\rho} |k_j(\mathbf{x},\bbeta)-
k_j(\mathbf{x}_i,\bbeta) |G(\ddr\bbeta)\\
&\leq&\gamma+2G\{\bbeta\dvt
|\bbeta|>M_\rho\}\leq\gamma+2(\rho+\delta).
\end{eqnarray*}
Similar arguments lead to $I_3\leq\gamma+2\rho$. Finally, we get
\[
\biggl|\frac{q_j(\mathbf{x};G)}{ q_j(\mathbf{x};G_0)}-1 \biggr|\leq\frac{3\delta+2\gamma
+4\rho}{ c},
\]
so that, for given $\eta<1$, it is always possible to choose
$\delta$, $\rho$ (by tightness of $G_0$) and $\gamma$ (by the
Arzel\`{a}--Ascoli theorem) small enough such that the right-hand
side in the last inequality is smaller than $\eta$. The proof is
then complete.
\end{pf}

We now aim to show that \eqref{eq:num} holds under the hypotheses
of Theorem \ref{theorem:consistency}, by extending the method set
forth by Walker (\citeyear{walk04}) for strong consistency. In order to simplify
the notation, let $\Lambda_{nj}=\int_{A_j}\Lambda_n(G)\pp(\ddr
G)$, where $(A_j)_{j\geq1}$ is the covering of $A$ in
\eqref{eq:cover}. The following identity is the key:
%
\begin{equation}\label{eq:key}
\Lambda_{n+1j}/\Lambda_{nj}=q^{nA_j}_{_{Y_{n+1}}}(\X_{n+1}) /
q_{_{Y_{n+1}}}(\X_{n+1};G_0),
\end{equation}
where $q^{nA_j}_{l}(\X_{n+1})=\int_{\mathbb{P}}
q_l(\X_{n+1};G)\pp_{nA_j}(\ddr G)$, $l\in\mathbf{C}$ and $\pp
_{nA_j}$ is
the posterior distribution restricted, and normalized, to the set
$A_j$. Note that \eqref{eq:key} includes the case of $n=0$ and
$\Lambda_{0j}=\pp(A_j)$.
By using conditional expectation, we have that
\begin{eqnarray*}
\mathrm{E} [\Lambda_{n+1j}^{1/2}\vert (Y_1,\X_1),\ldots,(Y_n,\X
_n),\X_{n+1} ]
&=&\Lambda_{nj}^{1/2}\sum_{l\in\mathbf{C}}\sqrt{q^{nA_j}_l(\X_{n+1})
q_l(\X_{n+1};G_0)}\\
&=&\Lambda_{nj}^{1/2} \bigl(1-h [\q^{nA_j}(\X_{n+1}),
\q(\X_{n+1};G_0) ] \bigr),
\end{eqnarray*}
where
$\q^{nA_j}(\X_{n+1})=[q^{nA_j}_1(\X_{n+1}),\ldots,q^{nA_j}_J(\X_{n+1})]$
and, for $\q_1,\q_2\in\Delta$,
\[
h(\q_1,\q_2)=1-\sum_{j\in\mathbf{C}}\sqrt{q_{1j}q_{2j}}.
\]
Note that $h(\q_1,\q_2)$ is a variation of the Hellinger distance
$\sqrt{\sum_{j\in\mathbf{C}}(q_{1j}^{_{1/2}}-q_{2j}^{_{1/2}})^2}$ on
$\Delta$ and that $h(\q_1,\q_2)\leq1$.
By taking the conditional expectation with respect to
$(Y_1,\X_1),\ldots,(Y_n,\X_n)$ only, we get the following identity:
%
\begin{equation}\label{eq:key2}
\textrm{E}\{\Lambda_{n+1j}^{1/2}\vert(Y_1,\X_1),\ldots,(Y_n,\X_n)\}
=\Lambda_{nj}^{1/2} \biggl(1
-\int_{\mathcal X}h [\q^{nA_j}(\mathbf{x}),
\q(\mathbf{x};G_0) ]M(\ddr\mathbf{x}) \biggr).
\end{equation}
Since the Hellinger distance and the Euclidean distance are
equivalent metrics in $\Delta$, it can be proven that, for
$(\q_n)_{n\geq1}\in\mathscr{Q}$ and $\q_0\in\mathscr{Q}$,
%
\begin{equation}\label{eq:equivL1H}
\int_{\mathcal X}h [\q_n(\mathbf{x}),
\q_0(\mathbf{x}) ]M(\ddr\mathbf{x})\to0\quad\mbox{if and only if}\quad
d(\q_n,\q_0)\to0.
\end{equation}
The equivalence in \eqref{eq:equivL1H} can be used to show that
$\int_{\mathcal X}h [\q^{nA_j}(\mathbf{x}),\q(\mathbf{x};G_0) ]M(\ddr\mathbf{x})$ is
bounded away from zero. In fact, take $G_j$ defined in
\eqref{eq:cover} and note that, by the triangle inequality,
\begin{eqnarray*}
\int_{\mathcal X}h [\q^{nA_j}(\mathbf{x}),\q(\mathbf{x};G_0)
]M(\ddr\mathbf{x})&\geq&
\int_{\mathcal X}h [\q(\mathbf{x};G_j),\q(\mathbf{x};G_0) ]M(\ddr\mathbf{x})\\
&&{}-\int_{\mathcal X}h [\q^{nA_j}(\mathbf{x}),\q(\mathbf{x};G_j) ]M(\ddr
\mathbf{x}).
\end{eqnarray*}
Since $d (\q(\cdot;G_j),\q(\cdot;G_0) )>\epsilon$,
\eqref{eq:equivL1H} ensures the existence of a positive constant,
say $\epsilon_2$, such that $\int_{\mathcal
X}h [\q(\mathbf{x};G_j),\q(\mathbf{x};G_0) ]M(\ddr\mathbf{x})>\epsilon_2$.
Now, choose
$\eta$ in \eqref{eq:cover} such that, for each $G\in A_j$,
$\int_{\mathcal
X}h [\q(\mathbf{x};G),\q(\mathbf{x};G_j) ]M(\ddr\mathbf{x})<\epsilon_2$,
where we have again
used \eqref{eq:equivL1H}. Since $\q^{nA_j}(\mathbf{x})$ does not
correspond exactly to a particular
$G\in A_j$, we use the convexity of the distance\vspace*{1pt}
$h [\q(\mathbf{x};G),\q(\mathbf{x};G_j) ]$ in its first argument to show that
$\int_{\mathcal X}h [\q^{nA_j}(\mathbf{x}),\q(\mathbf{x};
G_j) ]M(\ddr\mathbf{x})<\epsilon_2$. Note that, in fact, by Jensen's
inequality,
\begin{eqnarray*}
\int_{\mathcal X}h [\q^{nA_j}(\mathbf{x}),\q(\mathbf{x};G_j) ]M(\ddr\mathbf{x})
&=&\int_{\mathcal X} \biggl(1-\sum_{l\in\mathbf{C}}
\sqrt{\int_{\mathbb{P}} q_l(\X_{n+1};G)\pp_{nA_j}(\ddr G)
q_l(\mathbf{x};G_j)} \biggr)M(\ddr\mathbf{x})\\
&\leq&
\int_{\mathbb{P}}\int_{\mathcal X}h [\q(\mathbf{x};G),\q(\mathbf{x};G_j) ]
M(\ddr\mathbf{x}) \pp_{nA_j}(\ddr G)<\epsilon_2.
\end{eqnarray*}
Hence, there exists a $\epsilon_3>0$ such that $\int_{\mathcal
X}h [\q^{nA_j}(\mathbf{x}),\q(\mathbf{x};G_0) ]M(\ddr\mathbf{x})>\epsilon_3$.

From \eqref{eq:key2}, it now follows that
\[
E(\Lambda^{1/2}_{n+1j})<(1-\epsilon_3)^n\sqrt{\pp(A_j)}
\]
and an application of Markov's inequality leads to
\[
\mathrm{P} \biggl\{\sum_{j\geq1}\Lambda^{1/2}_{nj}>\exp(-nb) \biggr\}<\exp(nb)
(1-\epsilon_3)^n\sum_{j\geq1}\sqrt{\pp(A_j)}.
\]
Therefore, \eqref{eq:den} holds for any $b<-\log(1-\epsilon_3)$ from
an application of the Borel--Cantelli lemma, provided that the
following summability condition is satisfied:
%
\begin{equation}\label{eq:summability}
\sum_{j\geq1}\sqrt{\pp(A_j)}<+\infty.
\end{equation}
Lemma \ref{lemma:3} below shows that $\pp$ satisfies condition
\eqref{eq:summability} under the stated hypotheses and, in turn,
completes the proof of Theorem \ref{theorem:consistency}.

\begin{lemma}\label{lemma:3}
Let $H\in\mathbb{P}$ be the prior predictive distribution of $\pp$
and assume that condition \textup{(ii)} of Theorem \ref{theorem:consistency}
holds. Then \eqref{eq:summability} is verified.
\end{lemma}

\begin{pf}
The proof follows along the lines of arguments used by Lijoi, Pr\"{u}nster and
Walker (\citeyear{LijPruWal05}). Take $\delta$ to be any positive number in $(0,1)$
and $(a_n)_{n\geq1}$ any increasing sequence of positive numbers
such that $a_n\to+\infty$. Also, let $a_0=0$. Define
$C_n=\{\bbeta\dvt |\bbeta|\leq a_n\}$ and consider the family of
subsets of $\mathbb{P}$ defined by
%
\begin{equation}\label{eq:calB}
\mathbb B_{a_n,\delta}= \{G\dvt G(C_n)\geq1-\delta,
G(C_{n-1})<1-\delta\}
\end{equation}
for each $n\geq1$. These sets are pairwise disjoint and $\bigcup_n
\mathbb B_{a_n,\delta}=\mathbb{P}$. For the moment, let us assume
that the metric entropy of $\mathbb B_{a_n,\delta}$ with respect to
the distance $d$ is uniformly bounded in $n$, that is, the number of
$\eta$-balls in the distance $d$ that covers $\mathbb
B_{a_n,\delta}$ is finite for any $n$. Summability in
\eqref{eq:summability} is then implied by
%
\begin{equation}\label{eq:summability2}
\sum_{n\geq1}\sqrt{\pp(\mathbb B_{a_n,\delta})}<+\infty.
\end{equation}
In order to prove \eqref{eq:summability2},
note that $\mathbb B_{a_n,\delta}\subset\{G\dvt
G(C^c_{n-1})>\delta'\}$ for some $\delta'>\delta$.\vspace*{1pt} An application of
Markov's inequality leads to $\pp(\mathbb B_{a_n,\delta})\leq
(1/\delta')H(C^c_{n-1})$,
hence \eqref{eq:summability2} is implied by
$\sum_{n\geq1}\sqrt{H(C^c_{n-1})}<+\infty$. Next, we have that
\[
\int_{\rr^d}|\bbeta|H(\ddr\bbeta)=
\sum_{n\geq1}\int_{C^c_{n-1}/C^c_n}|\bbeta|H(\ddr\bbeta)
\geq\sum_{n\geq1}a_{n-1}[H(C^c_{n-1})-H(C^c_n)],
\]
by a second application of Markov's inequality, so that condition
(ii) of Theorem \ref{theorem:consistency} ensures that
$\sum_{n\geq1}a_{n-1}[H(C^c_{n-1})-H(C^c_n)]<+\infty$. If we now
take $a_n\sim n^2$, it is easy to see that $H(C^c_n)=o(n^{-(2+r)})$
for some $r>0$. For example,
\[
\sum_{n\geq1}(n-1)^2[H(C^c_{n-1})-H(C^c_n)]=\sum_{n\geq1} (2n-1)H(C^c_{n}).
\]
This, in turn, ensures the convergence of\vspace*{1pt}
$\sum_{n\geq1}H(C^c_{n-1})^\alpha$ for any $\alpha$ such that
$(2+r)^{-1}<\alpha<1$, which includes the case $\alpha=1/2$.
Condition \eqref{eq:summability2} is then verified.

In order to complete the proof, it remains to show that the metric
entropy of $\mathbb B_{a_n,\delta}$ with respect to the distance $d$
is uniformly bounded in $n$. It is actually sufficient to reason in
terms of the distance over $\mathbb{P}$ induced by
\[
d_j(\q_1,\q_2)=\int_{\mathcal X}|q_{1j}(\mathbf{x})-
q_{2j}(\mathbf{x})|M(\ddr\mathbf{x})
\]
for an arbitrary $j\in\mathbf{C}$ since $\max_{j}d_j(\q_1,\q
_2)\leq
d(\q_1,\q_2)\leq J\max_{j}d_j(\q_1,\q_2)$.
Let $\mathscr G$ be a set in $\mathscr{Q}$
and, for $\delta>0$, denote by $J(\delta,\mathscr G)$ the metric
entropy of $\mathscr G$ with respect to $d_j$, that is, the logarithm
of the minimum of all $k$ such that there exists
$\q_1,\ldots,\q_k\in\mathscr{Q}$ with the property that
$\forall\q\in\mathscr{Q}$, there exists an $i$ such that
$d_j(\q,\q_i)<\delta$.
The result is then stated as follows: for ${\mathscr
G}_{a_n,\delta}= \{\q(\mathbf{x};G)\dvt G\in\mathbb B_{a_n,\delta} \}$,
there exists an $M_\delta<+\infty$ depending only on $\delta$ such
that, for any $n$,
%
\begin{equation}\label{eq:entropy}
J(\delta,{\mathscr G}_{a_n,\delta})<M_\delta.
\end{equation}
The proof of \eqref{eq:entropy} consists of a sequence of three
steps.

%
%
\textit{Step (1)}. Define $C_a=\{\bbeta\dvt |\bbeta|\leq a\}$ and
$\mathscr F_a=\{\q(\mathbf{x};G)\dvt G(C_a)=1\}$. Then
%
\begin{equation}\label{eq:step1}
J(2\delta,\mathscr F_a)\leq \biggl(\frac{2aK}{\delta}+1 \biggr)^d
\biggl(1+\log\frac{1+\delta}{\delta} \biggr),
\end{equation}
where $K$ is a constant that depends on the total volume of the
space $\mathcal X$. It is easy to show that, for any $j\in\mathbf
{C}$, the
kernel $k_j(\mathbf{x},\bbeta)$ is a Lipschitz function in $\bbeta$ with
Lipschitz constant $K_{\mathbf{x}}=\max_{i\leq J}\{|\mathbf{x}_j-\mathbf{x}_i|\}$.
Hence,
\[
\int_{\mathcal X}|k_j(\mathbf{x},\bbeta_1)-
k_j(\mathbf{x},\bbeta_2)|M(\ddr\mathbf{x})\leq K|\bbeta_1-\bbeta_2|,
\]
where $K=\sup_{\mathbf{x}\in\mathcal X}K_{\mathbf{x}}<+\infty$. Given
$\delta$, let
$N$ be the smallest integer greater than $4aK/\delta$ and cover
$C_a$ with a set of balls $E_i$ of radius $2a/N$
so that, for any $\bbeta_1,\bbeta_2\in E_i$,
$|\bbeta_1-\bbeta_2|<4a/N$. This leads to $\int_{\mathcal
X}|k_j(\mathbf{x},\bbeta_1)-k_j(\mathbf{x},\bbeta_2)|M(\ddr\mathbf{x})\leq
\delta$. The
number of balls necessary to cover $C_a$ is then smaller than $N^d$.
Using arguments similar to those used in Ghosal, Ghosh and
Ramamoorthi (\citeyear{ghosal99}), Lemma 1, it can be shown that
$J(2\delta,\mathscr F_a)\leq
N^d (1+\log[(1+\delta)/\delta] )$, from which \eqref{eq:step1}
follows.

\textit{Step (2)}. Define $\mathscr F_{a,\delta}=\{\q(\mathbf{x};G)\dvt
G(C_a)\geq1-\delta\}$. Then
%
\begin{equation}\label{eq:step2}
J(\delta,\mathscr F_{a,\delta})\leq K_\delta a^d
\end{equation}
for a constant $K_\delta$ depending on $\delta$. To see this, take
$\q(\mathbf{x};G)\in\mathscr F_{a,\delta}$ and denote by $G^*$ the
probability measure in $\mathbb P$ defined by $G^*(A)=G(A\cap C_a)/
G(C_a)$ so that $\q(\mathbf{x};G^*)$ belongs to $\mathscr F_a$. It is easy
to verify that $d_j (\q(\cdot;G^*),\q(\cdot;G) )<2\delta$.
It follows that $J(3\delta,\mathscr F_{a,\delta})\leq
J(\delta,\mathscr F_a)$, from which \eqref{eq:step2} follows.

%
%
\textit{Step (3)}. We follow here a technique used by Lijoi,
Pr\"{u}nster and Walker (\citeyear{LijPruWal05}), Section 3.2. For the sequence
$(a_n)_{n\geq1}$ introduced before, define
\[
\mathscr F^U_{a_n,\delta}=\{\q(\mathbf{x};G)\dvt
G(C_n)\geq1-\delta\}\quad\mbox{and}\quad
\mathscr F^L_{a_n,\delta}=\{\q(\mathbf{x};G)\dvt G(C_n)<1-\delta\}.
\]
By construction, ${\mathscr G}_{a_n,\delta}\subset\mathscr
F^U_{a_n,\delta}$ and ${\mathscr G}_{a_n,\delta}\subset\mathscr
F^L_{a_{n-1},\delta}$. Moreover, $\mathscr
F^L_{a_{n-1},\delta}\downarrow\varnothing$ as $n$ increases to
$+\infty$, thus, for any $\eta>0$, there exists an integer $n_0$
such that, for any $n\geq n_0$, $J(\eta,\mathscr
F^L_{a_n,\delta})\leq J(\eta,\mathscr F^U_{a_{n_0},\delta})$. By
\eqref{eq:step2}, it follows that
%
\begin{equation}\label{eq:step3}
J(\eta,\mathscr G_{a_n,\delta})\leq K_\delta a_{n_0}^d
\end{equation}
for any $n\geq n_0$, but, since $\mathscr
G_{a_n,\delta}\subset\mathscr F^U_{a_n,\delta}$ and $\mathscr
F^U_{a_{n},\delta}\uparrow\mathscr{Q}$,
\eqref{eq:step3} is also true for any $n<n_0$. Result
\eqref{eq:entropy} is then verified by setting $M_\delta=K_\delta
a_{n_0}^d$.
\end{pf}


\section*{Acknowledgments}

The authors are grateful to an anonymous referee for
his valuable comments and suggestions which led to a substantial
improvement of the paper. Special thanks are also due to A.~Lijoi and
I.~Pr\"{u}nster for some useful discussions. P.~De Blasi was partially
supported by Regione Piemonte. J.W.~Lau's research was partly supported
by Hong Kong RGC Grant \#601707. L.F.~James was supported in part by
grants HIA05/06.BM03, RGC-HKUST 6159/02P, DAG04/05.BM56 and RGC-HKUST
600907 of the HKSAR.

\printhistory

\end{document}